\documentclass[a4paper]{article}

\usepackage{amsfonts}
\usepackage{latexsym}
\usepackage{amssymb}
\usepackage{amscd}

\renewcommand{\Sigma}{\sum}
\newtheorem{thm}{\sc theorem}[section]

\newtheorem{lemma}[thm]{\sc lemma}
\newtheorem{prop}[thm]{\sc proposition}
\newtheorem{cor}[thm]{\sc corollary}
\newtheorem{conj}[thm]{\sc conjecture}

\newenvironment{pf}{{\emph{Proof.}}}{\hfill$\Box$\\[1mm]}

\newcommand{\A}{{\mathcal{A}}}
\newcommand{\B}{{\mathcal{B}}}

 \def\b{\beta} \def\g{\gamma} \def\d{\delta}
\def\e{\varepsilon}   \def\l{\lambda}
  \def\s{\sigma}  \def\t{\tau}
\def\om{\omega}

\def\bC{\mathbb{C}}
\def\bN{\mathbb{N}}

\def\bZ{\mathbb{Z}}

\def\A{\mathcal{A}}

\def\sA{{}_\s \A}

\def\M{\mathcal{M}}
\def\N{\mathcal{N}}
\def\Eacute{\mathrm{\acute{E}}}
\def\eacute{\mathrm{\acute{e}}}

 \def\id{\mathrm{id}}

\def\der{\partial}
\def\im{\mathrm{im}}
\def\ker{\mathrm{ker}}

\def\Podles{\mathrm{Podle{\acute{s}}}}
\def\integers{\mathrm{integers}}

\def\linear{\mathrm{linear}}
\def\otherwise{\mathrm{otherwise}}

\def\even{\mathrm{even}}
\def\odd{\mathrm{odd}}

\def\b{\beta}

\def\suq2{SU_q (2)}
\def\csuq2{C(SU_q (2))}
\def\slq2{SL_q (2)}
\def\s{\sigma}
\def\smod{\sigma_{mod}}

\def\rto{\rightarrow}
\def\w2{\omega_2}
\def\w2p{\omega^{'}_2}

\def\cqg{\A(G_q)}

 \def\hhas{HH_\ast^\s(\A)}
\def\hcas{HC_\ast^\s(\A)}

\def\Asnzero{C^\s_0}
\def\Asnone{C^\s_1}
\def\Asntwo{C^\s_2}
\def\Asnthree{C^\s_3}
\def\Asnn{C^\s_n}
\def\Asnnplusone{C^\s_{n+1}}
\def\Asnstern{C^\s_\ast}

\begin{document}

\title{Twisted homology of quantum $SL(2)$}
\author{Tom~Hadfield\footnote{Supported until 31/12/2003 by the EU
Quantum Spaces - Noncommutative Geometry Network (INP-RTN-002) and
from 1/1/2004 by an EPSRC postdoctoral fellowship} ${}^1$,
Ulrich~Kr\"{a}hmer${}^2$} 
\date{\today}
\maketitle

\centerline{${}^1$ School of Mathematical Sciences,}
\centerline{Queen Mary, University of London}
\centerline{327 Mile End Road, London E1 4NS, England}
\centerline{t.hadfield@qmul.ac.uk}
\centerline{}
\centerline{ ${}^2$ Humboldt Universit\"at zu Berlin,}
\centerline{Institut f\"ur Mathematik}
\centerline{Unter den Linden 6}
\centerline{Sitz: Rudower Chaussee 25}
\centerline{D-10099 Berlin, Germany}
\centerline{kraehmer@mathematik.hu-berlin.de}
\centerline{}
\centerline{MSC (2000): 58B34, 19D55, 81R50, 46L}

\abstract{We calculate the twisted Hochschild and cyclic homology 
(in the sense of Kustermans, Murphy and Tuset) of the coordinate
algebra of the quantum $SL(2)$ group    
relative to twisting automorphisms acting by rescaling the
standard generators $a,b,c,d$. We discover a 
family of automorphisms for which the
``twisted" Hochschild dimension coincides with the classical
dimension of $SL(2, \bC)$, thus avoiding the ``dimension drop" in
Hochschild homology seen for many quantum deformations. 
Strikingly, the simplest such automorphism is the canonical modular
automorphism arising from the Haar functional.
In addition, we identify the twisted cyclic cohomology classes
corresponding to the three covariant differential calculi over quantum
$SU(2)$ discovered by Woronowicz.}

\section{Introduction}

Cyclic homology and cohomology were independently discovered 
by Alain Connes \cite{connes0} and Boris Tsygan \cite{tsygan} in the
early 1980's, and should be thought of as extensions of de Rham (co)homology 
to various categories of noncommutative algebras. Quantum groups also
appeared  in the same period, 
with the first example of a ``compact quantum group'' in the
C*-algebraic setting being Woronowicz's ``quantum $SU(2)$'' \cite{wo_suq2}.

The noncommutative differential geometry (in the sense of Connes) of 
quantum $SU(2)$ was thoroughly investigated by Masuda, Nakagami and 
Watanabe  \cite{mnw}. They calculated the Hochschild and cyclic 
homology of the coordinate algebra $\A(SL_q(2))$
of quantum $SL(2)$ as well as 
the K-theory and K-homology of the 
C*-algebra of the compact quantum $SU(2)$ group. 
This work was extended by Feng and Tsygan \cite{FT}, who
computed the Hochschild and cyclic homology of 
the standard quantized coordinate algebra $\cqg$ associated 
to an arbitrary complex semisimple Lie group $G$.
The  homologies are 
roughly speaking those  of a classical space labelling the symplectic
leaves of the Poisson-Lie group $G$ (the semi-classical
limit of $\cqg$). In particular, the
Hochschild dimension of $\cqg$
equals the rank of $G$. 
This ``dimension drop" had already been observed 
for other quantizations of  Poisson algebras.
Many authors regarded it as an unpleasant 
feature and asked for generalizations of cyclic homology 
which detect the quantized parts of quantum groups as well.

One candidate is  twisted Hochschild and cyclic (co)homology 
defined by Kustermans, Murphy and Tuset \cite{kmt}, relative to a pair of an algebra $\A$ and automorphism $\s$. This reduces to  ordinary Hochschild
and cyclic (co)homology of $\A$ on taking $\s$ to be the identity.
 The standard theory is intimately related with the
idea of considering tracial functionals on noncommutative algebras as
analogues of integrals, whereas the twisted theory arises naturally
from functionals whose tracial properties are of the form $h(ab)=h(\sigma(b)a)$. 
Noncommutative spaces equipped with such functionals include
 duals of nonunimodular groups, type III von Neumann algebras
and compact quantum groups. The aim of 
\cite{kmt}  was to adapt Connes' constructions
relating cyclic cohomology and differential calculi to
covariant differential calculi in the sense of Woronowicz, 
since the volume forms of such calculi 
define in general twisted cocycles rather than
usual ones \cite{sw1}. The possibility of 
pairing twisted cyclic cocycles (e.g.~over quantum homogeneous 
spaces) with equivariant K-theory was demonstrated in \cite{NT04}, and
it seems an interesting problem to investigate how far this original motivation of cyclic cohomology
extends to the twisted setting.

In this paper we compute the twisted Hochschild and cyclic homologies $\hhas$, $\hcas$
for the coordinate algebra $\A=\A(SL_q(2))$ of the quantum $SL(2)$ group, with generic deformation parameter $q$. 
We consider all automorphisms $\s$ of the form
$a,b,c,d \mapsto \lambda a, \mu b, \mu^{-1} c, \l^{-1} d$,
where $a, b, c, d$ are the standard generators, and $\l$, 
$\mu$ are nonzero elements of $k$. 
As an overview we 
collect the dimensions of $HH_n^\s(\A)$ as a $k$-vector space, 
see the main text for explicit formulas for generators:  

\begin{thm}\label{erga}
We have
\begin{eqnarray}
&& \mathrm{dim}\,HH_n^\s (\A)=0,\quad n > 3,\nonumber\\   
&& \mathrm{dim}\,HH_3^\s (\A)=\left\{
	\begin{array}{ll}
	N+1 \quad &\lambda=q^{-(N+2)},\mu=1,\\
	0 \quad &\otherwise,
	\end{array}\right.\nonumber\\
&& \mathrm{dim}\,HH_2^\s(\A)=\left\{
	\begin{array}{ll}
	N+1 \quad &
	\lambda=q^{-(N+2)},\mu=1,\\
	2 \quad & \lambda=q^{-(N+1)},\mu=q^{\pm(M+1)},\\
	0 \quad & \otherwise,
	\end{array}\right.\nonumber\\
&& \mathrm{dim}\,HH_1^\s (\A)=\left\{
	\begin{array}{ll}
	0 \quad & 
	\lambda \notin q^{-\bN},
	\mu=q^{\pm(M+1)},\\
	0 \quad & 
	\lambda \neq 1,\mu \notin q^\mathbb{Z},\\
	4 \quad &
	\lambda=q^{-(N+1)},\mu=q^{\pm(M+1)},\\
	\infty \quad & \otherwise,
	\end{array}\right.\nonumber\\
&& \mathrm{dim}\,HH_0^\s (\A)=\left\{
	\begin{array}{ll}
	\infty \quad & \mu=1,\\
	2 \quad & \lambda=q^{-(N+1)},\mu=q^{\pm(M+1)},\\
	0 \quad & \otherwise,
	\end{array}\right.\nonumber 
\end{eqnarray}
for $M, N \in \bN$.
\end{thm} 

Strikingly (Theorem \ref{HH_3_for_any_sigma}), there exists a  family of automorphisms for 
which the twisted Hochschild dimension takes the classical
value three (note also that the homological dimension of $\A(SL_q(2))$ is
three \cite{ls}) -  the twisted theory avoids the ``dimension drop". Remarkably,
the simplest such automorphism  
($\l = q^{-2},\mu=1$) 
is  the canonical modular automorphism associated to the 
Haar functional on $\A$. Similar results were obtained for $\Podles$
quantum spheres \cite{tom} and quantum hyperplanes \cite{si}.

In \cite{FT}, Feng and Tsygan considered formal quantizations,  with
 $\cqg$ a Hopf algebra over 
$\bC [[\hbar]]$ with $q=e^\hbar$. They showed 
that for a Hopf algebra $\A$ over a field $k$, with coproduct $\Delta$, counit 
$\varepsilon$ and antipode $S$, and an $\A$-bimodule $\M$, 
there is an isomorphism 
\begin{equation}
H_n(\A,\M) \simeq 
\mathrm{Tor}_n^\A(\M',k)
\end{equation}
 Here,
$\M'$ is $\M$ as a linear space with right
action given by 
\begin{equation}
m \blacktriangleleft a:= \sum \; S(a_{(2)})ma_{(1)}
\end{equation}
using Sweedler's notation for the coproduct, 
and $k=\A/ \ker \,\varepsilon$ is the 
trivial left $\A$-module. Then they computed these $\mathrm{Tor}$-groups
using the spectral sequence associated to the filtration
induced by $\hbar$.

In this paper we  compute
$\mathrm{Tor}_n^\A(\M',k)$
from a Koszul-type free resolution
\begin{equation}
	0 \rightarrow \A \rightarrow \A^3 \rightarrow \A^3 
	\rightarrow \A \rightarrow k \rto 0
\end{equation}
of $k$. Noncommutative Koszul resolutions were studied
by several authors, in particular Wambst \cite{wambst}, but as far
as we know were not applied to quantum groups.
 In our opinion this resolution shows very clearly the 
geometric mechanisms behind the computations. 
We will see that the maps of the resulting complex computing the 
twisted Hochschild homology become zero for $q=1$, so  
one obtains the Hochschild-Kostant-Rosenberg theorem for
$SL(2)$ (the algebraic cotangent bundle of $SL(2)$ is
trivial). However, for $q \neq 1$ this does not happen 
for any twisting automorphism.

A summary of this paper is as follows.
In section  \ref{section_thch} we  recall how the twisted theory was discovered \cite{kmt}, then give the definitions of $HH_\ast^\s(\A)$ and 
$HC_\ast^\s(\A)$, and the underlying cyclic object.  
We specialize to Hopf algebras and 
explain the methods adapted from \cite{FT}. We then
present the general scheme of the
noncommutative Koszul complexes used here. 
In section \ref{section_slq2_suq2} we introduce the   
quantum $SL(2)$ group.
In section \ref{section_th_slq2} we present our 
calculations of $HH_\ast^\s(\A)$ for $\A=\A(SL_q (2))$.

Twisted cyclic homology is defined as the total homology of Connes'
mixed $(b,B)$-bicomplex coming from the underlying cyclic object, as
in \cite{loday}. In section \ref{section_tc_slq2} we compute this
homology via a spectral sequence.

Finally, in section~\ref{section_diff_calc_suq2} we discuss 
the relation of our results to previously known twisted 
cyclic cocycles coming from the three covariant 
differential calculi over $\A(SL_q(2))$ discovered by Woronowicz.
The twisted cyclic 3-cocycle arising from the three dimensional left 
covariant calculus was given explicitly in \cite{kmt} and \cite{sw1}. 
We show (Theorem \ref{3dcalc_cocycle_is_trivial}) that this 
3-cocycle is a trivial element of twisted cyclic cohomology.
Further, the twisted 4-cocycles arising from the two bicovariant four 
dimensional calculi both correspond to the twisted 0-cocycle coming 
from the Haar functional (as elements of even periodic twisted cyclic cohomology).

\section{Twisted cyclic homology}\label{section_thch}

\subsection{\sc motivation}
Twisted cyclic (co)homology arose from the study of covariant 
differential calculi over quantum groups \cite{kmt}. 

Let $\A$ be an algebra over $\bC$. 
Given a differential calculus
$(\Omega, d)$ over $\A$,  with 
$\Omega = \oplus_{n=0}^N \; \Omega_n$, Connes \cite{connes1} considered linear
functionals $\int : \Omega_N \rto \bC$, which are closed and graded traces on
$\Omega$, meaning 
$$\int \; d \omega = 0 \quad \forall \; \omega \in \Omega_{N-1}$$
\begin{equation}
\label{closed}
\label{graded_trace}
\int \; \omega_m \omega_n = (-1)^{mn} \int \; \omega_n \omega_m \quad \forall \; \omega_m \in \Omega_m, \; \omega_n \in \Omega_n
\end{equation}
Connes found that such linear functionals are in one to one
correspondence with cyclic $N$-cocycles $\tau$ on the algebra, via 
\begin{equation}\label{cyclic_cocycle}
	\tau( a_0, a_1, \ldots , a_N) = \int \; a_0 \; d a_1 \; d a_2 \ldots d a_N
\end{equation}
which led directly to his simplest formulation of cyclic cohomology \cite{connes1}.

If $\A$ is the coordinate algebra of a quantum group, then 
Woronowicz proposed to study covariant differential calculi, for which the 
left coaction of $\A$ on $\A$ given by the coproduct 
$\Delta : \A \rightarrow \A \otimes \A$ 
extends to a coaction $\Delta_L : \Omega \rto \A \otimes \Omega$
compatible with the differential $d$ \cite{wo_suq2}, \cite{wo_cmpg}. 
For such calculi the natural linear functionals
$\int : \Omega_N \rto \bC$ are 
no longer graded traces, but twisted graded traces, meaning that 
\begin{equation}
\label{twisted_graded_trace}
\int \; \omega_m \omega_n = (-1)^{mn} \int \; \s(\omega_n) \omega_m \quad \forall \; \omega_m \in \Omega_m, \; \omega_n \in \Omega_n
\end{equation}
for some degree zero automorphism $\s$ of $\Omega$.
In particular, $\s$ restricts to an automorphism of $\A$, and, for any $a \in \A$, $\omega_N \in \Omega_N$ we have
\begin{equation}
\label{int_aut}
\int \; \omega_N a = \int \; \s(a) \omega_N
\end{equation}
Hence for each covariant calculus there is a natural automorphism of $\A$.
Motivated by this observation, Kustermans, Murphy and Tuset defined
``twisted'' Hochschild and cyclic cohomology for any pair 
of an algebra $\A$ and automorphism $\s$, and showed that 
the one-to-one correspondence between graded traces and cyclic cocycles
generalizes to this setting. The next section recalls their 
definitions, transposed to homology.

\subsection{\sc twisted hochschild and cyclic homology} 
\label{section:cyclic_module}

Let $\A$ be a unital, associative algebra over a field $k$ (assumed to be of characteristic zero) 
and $\sigma$ an automorphism. We  define the cyclic object \cite{connes2}, \cite{loday} underlying
twisted cyclic homology $HC_\ast^\s (\A)$ 
of $\A$ relative to $\s$. Set $C_n:=\A^{\otimes (n+1)}$.
For clarity, we will denote
$a_0 \otimes a_1 \otimes \cdots \otimes a_n \in C_n$ 
by $(a_0, a_1, \ldots , a_n)$. Define
$$d_{n,i} (a_0, a_1, \ldots , a_n) = (a_0, \ldots , a_i a_{i+1} , \ldots ,
a_n) \quad 0 \leq i \leq n-1 $$ 
$$
	d_{n,n} (a_0, a_1, \ldots , a_n) = (\s(a_n) a_0, a_1, \ldots , a_{n-1})
$$
$$
	s_{n,i} (a_0, a_1, \ldots , a_n) = 
	(a_0 , \ldots a_i , 1 , a_{i+1} , \ldots , a_n) \quad 0 \leq i \leq n
$$
\begin{equation}
\t_n (a_0, a_1, \ldots , a_n) = ( \s(a_n) , a_0 , \ldots , a_{n-1})
\end{equation}
For $ \sigma = \id$ these are 
the face, degeneracy and cyclic operators of the
standard cyclic object associated to $\A$  \cite{loday}. 
For general $ \sigma $ the operator $ T_n:=\t_n^{n+1}$
is not equal to the identity, but all other relations of the cyclic
category are fulfilled. Hence $C_\ast$ becomes 
what is called a paracyclic object \cite{getjon}.
To obtain a cyclic object, we  pass to the cokernels 
$\Asnn := C_n / C_n^1$, $C_n^1:=\im (\id-T_n)$.
Dualizing \cite{kmt}, we call the cyclic homology of this cyclic
object the $\s$-twisted cyclic homology $HC_\ast^\s (\A)$ of $\A$.  
Hence $HC_\ast^\s (\A)$ is the total homology of Connes' mixed 
$(b,B)$-bicomplex

\begin{equation}\label{gemischt}
\label{mixed_b_B_bicomplex}
     	\begin{CD}
	@ V{b_4} VV @ V{b_3} VV @ V{b_2} VV @ V{b_1} VV @ . @ . @ .\\
	{\Asnthree} @ <{B_2}<< {\Asntwo} @ <{B_1} << {\Asnone} @ <{B_0} <<
	{\Asnzero} @ . @ . @ . @ .\\ 
 	@ V{b_3} VV @ V{b_2} VV @ V{b_1} VV @ . @ . @ . @ .\\
	{\Asntwo} @ <{B_1} << {\Asnone} @ <{B_0} << {\Asnzero} @ . @ . @ . @ . @ .\\
 	@ V{b_2} VV @ V{b_1} VV @ . @ . @ . @ . @ .\\
	{\Asnone} @ <{B_0} << {\Asnzero} @ . @ . @ . @ . @ . \\
 	@ V{b_1} VV @ . @ .  @ . @ . @ . @ .\\
	{\Asnzero} @ . @ . @ . @ . @ . @ .\\
    	\end{CD}
\end{equation}\\

The maps $b_n$ and $B_n$ are given by 
\begin{equation}
b_n = \Sigma_{i=0}^n \; (-1)^i d_{n,i},\quad
	B_n = (1+ (-1)^n \t_{n+1} ) s_n N_n,
\end{equation}
with 
$N_n = \Sigma_{j=0}^n \; (-1)^{nj} \t_n^j$, and 
$s_n  :  \Asnn \rto \Asnnplusone$ the ``extra degeneracy''
\begin{equation}
s_n ( a_0, a_1, \ldots , a_n) = (1, a_0, a_1, \ldots , a_n)
\end{equation}

We calculate $HC_\ast^\s(\A)$ via the spectral sequence associated to 
the mixed complex. Let $HH_\ast^\s(\A)$
denote the entries of its first page, that is, 
$HH_n^\s(\A):=H_n(\Asnstern,b_\ast)$ (the homologies of the columns). 
For $\s=\id$
these are the Hochschild homologies $HH_\ast ( \A)=H_\ast ( \A,\A)$. 
Hence we call $HH_\ast^\s (\A)$ as
in \cite{kmt} the
$ \sigma $-twisted Hochschild homology of $\A$.

To compute $HH_\ast^\s (\A)$ consider the mixed complex (\ref{gemischt}) with $\Asnn$ replaced by the
original $C_n$. This is not a bicomplex: the commutation relations in a
paracyclic object imply that the
(lifts of the) operators $b_\ast$ and $B_\ast$ anticommute according
to (see \cite{getjon}, Theorem~{2.3})
\begin{equation}\label{neuegl}
	b_{n+1} B_n +B_{n-1} b_n=\id -T_n.
\end{equation}  
But the columns form the complex $(C_\ast,b_\ast)$ which
computes the Hochschild homology $H_\ast(\A,\sA)$ of $\A$ 
with coefficients in the bimodule $\sA$ which is 
$\A$ as a vector space with bimodule structure
\begin{equation}
\label{twisted_bimodule}
a \triangleright b  \triangleleft c:=\sigma(a)bc
\end{equation}
In many cases   
$C_n=C_n^0 \oplus C_n^1$, $C_n^0:=\mathrm{ker}(\id -T_n)$,
for example when $ \sigma $ is diagonalizable. In this
case, $T_n = \s^{\otimes (n+1)}$ 
is also diagonalizable, and 
$C_n^0$ and $C_n^1$ are the eigenspace of 
$T_n$ corresponding to the eigenvalue $1$ and the direct sum
of all other eigenspaces, respectively. Then:

\begin{prop}
\label{tatata}
If $C_n=C_n^0 \oplus C_n^1$, then
$H_\ast(\A,\sA) \cong HH_\ast^\s(\A)$.
\end{prop} 
\begin{pf}
Note that (\ref{neuegl}) implies that 
$b_\ast$ commutes with $\id -T_\ast$, so the
decomposition $C_n=C_n^0 \oplus C_n^1$ defines
a decomposition of complexes, and we can identify 
$HH_\ast^\s(\A)$ with the homologies of the subcomplex 
$(C_\ast^0,b_\ast) \subset (C_\ast,b_\ast)$.
Hence $H_\ast(\A,\sA)$ is the direct sum of $HH_\ast^\s(\A)$
and the homologies of $(C_\ast^1,b_\ast)$. But  
$ (\id - T_n)|_{C_n^1}$ is a bijection under these assumptions,
and we have on $C_n^1$ again by (\ref{neuegl}) the relation
$$
	b_{n+1} (1-T_n)^{-1} B_n +(1-T_n)^{-1} B_{n-1} b_n=\id.
$$
So $(\id - T_n)^{-1}B_\ast$ is a 
contracting homotopy for $(C_\ast^1,b_\ast)$ and the claim follows.
\end{pf}

This will allow us to calculate $HH_\ast^\s(\A)$ using standard
techniques of homological algebra.

The spectral sequence calculation is most efficiently done by passing to
the normalized mixed complex 
(see for example \cite{weibel}, Application~{9.8.4}). This leaves the
first page unchanged. The second step is to calculate the
horizontal homology of the rows relative to the maps $B_n$ which in the
normalized complex are given explicitly by 
\begin{equation}\label{B_n}
		B_n(a_0,\ldots,a_n) =
		\sum_{i=0}^{n} (-1)^{ni}
		(1,\sigma (a_i),\ldots,\sigma (a_n),a_0,\ldots,a_{i-1}).   
\end{equation} 
For quantum $SL(2)$, we find that everything stabilises at the second page,
and we can then read off the twisted cyclic homology.

For later use we note that by using the Hochschild-Kostant-Rosenberg
theorem applied to an appropriate subalgebra, we obtain:
\begin{lemma}
\label{lemma:B_0_x^s_y^t}
If $x$, $y$ are commuting elements of $\A$, with $\s(x) =x$, 
$\s(y)=y$, then for any $s$, $t \geq 0$ we have
$$
	B_0 [ x^s y^t ] = t [( x^s y^{t-1}, y)] + s[( x^{s-1} y^t ,x)] \in 
	HH_1^\s (\A)
$$
\end{lemma}

From now on, we will drop the suffices and write $b_n$  as $b$.

\subsection{\sc hochschild homology of hopf algebras} 

For arbitrary algebras, the Hochschild homologies are derived functors
in the category of $\A \otimes \A^\mathrm{op}$-modules, and working
with explicit resolutions usually involves lengthy calculations.
But if $\A$ is a Hopf algebra 
 then we can
describe $H_\ast(\A,\M)$ for an arbitrary $\A$-bimodule $\M$ as a derived
functor in the category of $\A$-modules.
Define a right $\A$-module $\M'$ which is $\M$ as a vector space
with right action given by 
\begin{equation}\label{act}
	m \blacktriangleleft a :=  \sum \; S(a_{(2)}) \triangleright m \triangleleft
	a_{(1)},\quad
	a \in \A,m \in \M.
\end{equation}  
Consider $k$ as the trivial $\A$-module 
$\A/\mathrm{ker}\, \varepsilon$. Feng and Tsygan proved:
\begin{prop}
\cite{FT}
\label{tor}
There is an isomorphism of vector spaces
$$
	H_n(\A,\M) \simeq 
	\mathrm{Tor}_n^\A(\M',k).
$$
\end{prop}
\begin{pf}
The $\mathrm{Tor}_n^\A(\M',k)$
are 
computed from the complex $(C_\ast,d)$ (with zeroth tensor component
now being $\M'$) with boundary map $d$ given by
\begin{equation}
d = \tilde d_0 + \sum_{i=1}^{n-1} (-1)^i d_i + (-1)^n \tilde d_n,
\end{equation}
where the $d_i$ are as above and 

$$	\tilde d_0 (a_0, a_1,  \ldots , a_n) 
:= (a_0 \blacktriangleleft a_1,  a_2,  \ldots , a_n),$$	
\begin{equation}
	\tilde d_n (a_0, a_1,  \ldots , a_n) 
:= (\varepsilon(a_n)a_0 , a_1 , \ldots , a_{n-1}).
\end{equation}
We define two linear maps $\xi,\xi' : C_n \rightarrow C_n$ by
$$ 
	\xi(a_0, a_1, \ldots , a_n)  
	:=( S((a_1  \ldots a_n)_{(2)}) \triangleright a_0,   
	(a_1)_{(1)}, \ldots , (a_n)_{(1)})$$
\begin{equation}\label{fetsimo}
\xi'(a_0,  \ldots , a_n)
:= ((a_1  \ldots  a_n)_{(2)}) \triangleright a_0,  
	(a_1)_{(1)}, \ldots , (a_n)_{(1)}).
\end{equation}
Then 
$\xi \circ \xi'=\xi' \circ \xi=\mathrm{id}_{C_n}$. It is easily
checked that $\xi$ commutes with $d_i$ for $1 \le i \le n-1$ and
that $\xi \circ \tilde d_i=d_i \circ \xi$, $i=0,n$. Hence 
$\xi \circ d=b \circ \xi$ and
$\xi$ is an isomorphism of complexes of $k$-vector spaces.
\end{pf}

Let $ \pi : \M' \rightarrow H_0(\A,\M) $ be the canonical 
projection. Then we have 
$ \pi (m \blacktriangleleft a)=\varepsilon(a) \pi(m)$, and
if we consider $H_0(\A,\M)$ as trivial right $\A$-module, then 
$ \pi \otimes \id_{\A^{\otimes n}}$ induces a morphism
$$
	H_n(\A,\M) 
 	\rightarrow 
	H_0(\A,\M) \otimes_k \mathrm{Tor}_n^\A(k,k). 
$$
If $\A$ is commutative and 
$\M=\A$ with the standard bimodule structure, then 
$H_0(\A,\M)=\A$, $\pi$ is the identity, and the above map is
the isomorphism of the Hochschild-Kostant-Rosenberg theorem. 
For  $\M=\sA$ the map defines a ``classical shadow" of 
twisted Hochschild homology.

\subsection{\sc noncommutative koszul resolutions}
\label{koszull}

Propositions~\ref{tatata} and~\ref{tor} allow us
to compute $HH_\ast^\s(\A)$ for Hopf algebras $\A$ 
and diagonalizable $ \sigma $ from a resolution of the trivial $\A$-module 
$k$. In the commutative case, 
such a resolution can be constructed in form
of a Koszul complex associated to a minimal set of generators of 
$\ker \, \e$. 
We will see that one can proceed in the same way 
for quantum $SL(2)$.
The general scheme of the construction of the resolution is as follows.

Let $\A$ be an algebra and $d$ be a positive integer.
Let $x_{i,j}$, $1 \leq i, j \leq d$ be elements of $\A$
satisfying
\begin{equation}\label{braid}
	x_{i,j}x_{i-1,k}=x_{i,k}x_{i-1,j}.%\quad j>k	
\end{equation} 
In the commutative case one can take  $x_{i,j}=x_{1,j}$, and in many examples 
the $x_{i,j}$ will be uniquely determined by the $x_{1,j}$.

For $0 \leq n \leq d$, let $K_n(x_{i,j})$ be the $\A$-bimodule
$\A^{(\!\!\!\tiny{\begin{array}{c} d \\ n\end{array}}\!\!\!)}$, which we 
identify for $n>0$ with the submodule of 
$\A^d \otimes_\A \cdots \otimes_\A \A^d$ ($n$ factors) 
spanned over $\A$ by
$$
	e_{i_1} \otimes_\A \cdots \otimes_\A e_{i_n},\quad
	1 \le i_1 < \ldots < i_n \le d,
$$
where $e_i$ is a basis of $\A^d$. For $n>d$ we set $K_n(x_{i,j}):=0$.
%Note that  $K_0 ( x_{i,j}) = \A$.\\
For an $\A$-bimodule $\N$, set 
$K_n(x_{i,j},\N):=
K_n(x_{i,j}) \otimes_\A \N$ and 
define $\A$-module maps 
$$
	k_m : K_n(x_{i,j},\N) \rightarrow
	K_{n-1}(x_{i,j},\N),\quad m=1,\ldots,n
$$ 
(we suppress the index $n$ at the $k_m$) by
$$
	k_m : e_{i_1} \otimes_\A \cdots \otimes_\A e_{i_n} \otimes_\A y \mapsto
	e_{i_1} \otimes_\A \cdots
	 \hat e_{i_m}  \cdots \otimes_\A e_{i_n} \otimes_\A
	y \triangleleft x_{n,{i_m}}.
$$
Then  for $r<s$:
\begin{eqnarray}
&&(k_r k_s - k_{s-1} k_r)
	(e_{i_1} \otimes_\A \cdots \otimes_\A e_{i_n} \otimes_\A y) \nonumber\\
&=& e_{i_1} \otimes_\A \cdots  \hat e_{i_r}  \cdots
	 \hat e_{i_s}  \cdots \otimes_\A e_{i_n} \otimes_\A 
	y \triangleleft (x_{n,i_s}x_{n-1,i_r}-x_{n,i_r}x_{n-1,i_s}).\nonumber
\end{eqnarray}
The last bracket vanishes by (\ref{braid}). Thus we get
\begin{prop}\label{serrek}
The map $k:=\sum_{r=1}^{n} (-1)^r k_r$ makes
$K_\ast(x_{i,j},\N)$ into a complex which we  
call the Koszul complex associated to $x_{i,j}$ and $\N$. 
\end{prop}

The zeroth homology of this complex is obviously 
the quotient of $\N$ by the submodule generated by all 
elements of the form $y \triangleleft x_{1,i}$, $y \in \N$, 
$i=1,\ldots,d$. The classical application of Koszul complexes is to
produce resolutions of this quotient, but the Koszul complex is not
always acyclic (see \cite{serre} for the commutative case).  
In our application we will take $\N=\A$ to be a Hopf algebra 
with the standard bimodule structure, and the 
$x_{1,j}$ ($1 \leq j \leq d$) will generate 
$\mathrm{ker}\,\varepsilon$ as a an (left or right) $\A$-module. The 
associated Koszul complex will be checked by hand to be 
acyclic (see Proposition~\ref{resolution} below), 
so it provides a resolution of $\A/\mathrm{ker}\, \varepsilon$, 
and $\mathrm{Tor}^\A_\ast(\M',k)$ equals the homologies of
the complex $\M' \otimes_\A K_\ast(x_{i,j},\A)$.
The quasi-isomorphism from this complex 
to the standard complex $(C_\ast,d)$ calculating the 
$\mathrm{Tor}$-groups described in the proof 
of Proposition~\ref{tor} is then given by
\begin{equation}
\label{morphism_of_complexes}
	e_{i_1} \otimes_\A \cdots \otimes_\A e_{i_n} \mapsto
	x_{1,i_1} \wedge \ldots \wedge x_{1,i_n} :=
	\sum_{s \in S_n} (-1)^{|s|}
	x_{n,i_{s(n)}} \otimes \cdots \otimes x_{1,i_{s(1)}}.  
\end{equation}

\section{Quantum $SL(2)$}
\label{section_slq2_suq2}

In this section, we introduce the main facts on the standard quantized
coordinate ring $\A=\A (SL_q (2))$ that will be used below. 

\subsection{\sc the hopf algebra $\A (SL_q (2))$}
\label{section_sl2su2}

Let $k$ be a field of characteristic zero, and $q \in k$ some nonzero
parameter, which we assume is not a root of unity. 
The coordinate algebra $\A=\A(\slq2)$ of the quantum group
$\slq2$ over $k$ is the $k$-algebra generated by symbols $a$, $b$,
$c$, $d$ with relations
$$ab=qba, \quad
ac=qca, \quad
bd=qdb, \quad
cd=qdc, \quad
bc=cb$$
\begin{equation}
\label{ad=}
ad-qbc=1, \quad
da - q^{-1} bc=1
\end{equation} 
There is a unique Hopf algebra structure on $\A$ such that 
$$ \Delta(a)=a \otimes a + b \otimes c,\quad \Delta(b)=a \otimes b+b \otimes d,$$
$$\Delta(c)=c \otimes a+d \otimes c,\quad \Delta(d)=c \otimes b+d \otimes d,$$
$$\varepsilon(a)=\varepsilon(d)=1,\quad \varepsilon(b)=\varepsilon(c)=0,$$
\begin{equation}
S(a)=d,\quad S(b)=-q^{-1}b,\quad S(c)=-qc,\quad S(d)=a.
\end{equation}
A vector space basis of $\A$ is given by the monomials
\begin{equation}\label{ZNN_filtration}
	e_{i,j,k}:=a^i b^j c^k,\quad i \in \bZ, \; j,k \in \bN,\quad
	a^i:=d^{-i}\quad\mbox{for}\quad i<0,	
\end{equation}
(we use the convention that $x^0 =1$, for $x \in \A$, $x \neq 0$). 
We have
$$
	e_{i,j,k}e_{l,m,n}=q^{-l(j+k)}e_{i+l,j+m,k+n}+
	\sum_{r>0} \lambda_{i,j,k,l,m,n}(r) e_{i+l,j+m+r,k+n+r}  
$$ 
for some constants $\lambda_{i,j,k,l,m,n}(r)$. It follows that
$\A$ admits a $\bZ$-grading and three separating 
decreasing $\bN$-filtrations
\begin{equation}
\A=\bigoplus_{i \in \bZ} \A^a_i,\quad
\A=\A^x_0 \supset \A^x_1 \supset \ldots,\quad
x=b,\; c,\; bc,
\end{equation} 
where $\A^a_i=\mathrm{span}\{e_{i,j,k}\}_{j,k}$ 
and $\A^x_n$ is the span of $e_{i,j,k}$ with $j,k,j+k \ge n$ for 
$x=b,c,bc$, respectively. For  $x \in \A$, let 
$x_i$ be its component in $\A^a_i$. Set $\A^x_{i,n}:=\A^x_n \cap
\A^a_i$. Then $\A^x_{i,n}\A^x_{j,m}=\A^x_{i+j,n+m}$. 
Define a Hermitian
inner product on $\A$ by requiring that $e_{i,j,k}$ are orthonormal and
 let $\pi_x$, $\pi_{i,j,k}$,
denote the orthogonal projections onto $(\A^x_1)^\perp$, $e_{i,j,k}$. We freely consider
$\pi_{i,j,k}$ as a map $\A \rightarrow k$.  
Note that $\pi_x(y)_i=\pi_x(y_i)$ for all $y \in \A$.

Finally, $\A$ has a $\bZ^2$-grading given by
\begin{equation}
\label{Z^2_grading}
\deg( e_{i,j,k}) = (i, j-k)
\end{equation}
This grading extends to $\A^{\otimes (n+1)}$ and is preserved by the
Hochschild boundary and the maps $B_n$ (\ref{B_n}). Hence $HH_\ast^\s (\A)$ and $HC_\ast^\s (\A)$ are
naturally $\bZ^2$-graded. 

\subsection{\sc the haar functional}
The Hopf algebra $\A$ is cosemisimple \cite{KS1}, that is, there is a 
unique linear functional $h : \A \rightarrow k$ 
satisfying $h(1)=1$ and
\begin{equation}
	(h \otimes \id ) \Delta (x) = h(x) 1 = ( \id \otimes h ) \Delta(x) 
	\quad \forall \; x \in \A
\end{equation} 
If $k=\mathbb{C} $ and $q \in \mathbb{R} $, then 
$\A$ can be made into a Hopf $\ast$-algebra whose C*-algebraic completion is Woronowicz's quantum $SU(2)$ group \cite{wo_suq2}. 
The functional $h$ extends to the Haar state of this compact quantum group. Hence (with slight abuse of terminology) we also in the general case call $h$ the
Haar functional of $\A$.
For any $x$, $y \in \A$, we have 
\begin{equation}
\label{defn_mod_aut}
h( xy ) = h ( y \smod(x))
\end{equation}
where $\smod$ is the so-called modular automorphism of $\A$ given by
\begin{equation}
\label{smod_slq2}
\smod (a) = q^{-2} a, \quad
\smod (d) = q^{2} d, \quad
\smod (b) = b, \quad
\smod (c) = c
\end{equation}
So $h$ is a 
$\smod^{-1}$-twisted cyclic 0-cocycle.
 
\subsection{\sc the automorphism group of $\A( SL_q (2))$} 
For $\l,\mu \in k^\times$ there are
unique automorphisms $\s_{\l, \mu},\t_{\l, \mu}$ 
of $\A$ with
$$ \sigma_{\lambda, \mu}(a)=\lambda a,\quad
	\sigma_{\lambda, \mu}(b)=\mu b,\quad
	\sigma_{\lambda, \mu}(c)=\mu^{-1} c,\quad
	\sigma_{\lambda, \mu}(d)=\lambda^{-1} d,$$
\begin{equation}
\label{defn_lambda_tau}
 \tau_{\lambda, \mu}(a)=\lambda a,\quad
	\tau_{\lambda, \mu}(b)=\mu^{-1} c,\quad
	\tau_{\lambda, \mu}(c)=\mu b,\quad
	\tau_{\lambda, \mu}(d)=\lambda^{-1} d.
\end{equation} 
It is easy to check that this list is complete, 
although we do not know a reference where this was
pointed out explicitly:
\begin{prop}\label{autos}
If $\sigma$ is an automorphism of 
$\A (SL_q (2))$, then either 
$\sigma=\sigma_{\lambda, \mu}$ or
$\sigma=\tau_{\lambda, \mu}$ for some $\lambda,\mu$.
\end{prop}  
\begin{pf}
Using the $\bZ$-grading and the 
$\bN$-filtrations mentioned above it is a straightforward calculation
to check that up to rescaling and exchanging $b$ and $c$ the original
generators are the only elements of the algebra that fulfill the
defining relations. 
\end{pf}

The  $\s_{\l,\mu}$ act diagonally with respect to the generators $a,b,c,d$.
The $\t_{\l,\mu}$ are also diagonalizable.
For fixed $\l, \mu$ define 
$ x_\pm = c \pm \mu b$.
Then $\t_{\l,\mu} ( x_{\pm} ) = \pm x_{\pm}$, and $a, x_+ , x_- ,d$  generate $\A$.
So by Proposition \ref{tatata}:

\begin{cor}
\label{crucial_isomorphism}
For $\A = \A( SL_q (2))$, and for each $n$ and every automorphism $\s$, we have
$HH_n^\s (\A) \cong H_n (\A, \sA)$. 
\end{cor}

\section{Twisted Hochschild homology of $\A( SL_q (2))$}
\label{section_th_slq2}

\subsection{\sc a koszul resolution of $\A/ \ker \, \varepsilon$}
Using the above facts it is easy to see that $\ker \,\e$ 
is generated as both a left and right $\A$-module by 
$x_{1,1}:=a-1,\; x_{1,2}:=b, \;x_{1,3}:=c$. 
For these elements there exists a Koszul complex 
$(K_\ast,k)$, $K_n:=K_n(x_{i,j},\A)$, in the 
sense of section~\ref{koszull} with the $x_{i,j}$ given by
\begin{equation}
	\left(\begin{array}{ccc} 
	a-1 & b & c\\
	q^{-1}a-1 & b & c\\
	q^{-2}a-1 & b & c 
	\end{array}\right).  
\end{equation}

We check by explicit calculation that this Koszul complex
is acyclic:
\begin{prop}\label{resolution}
The left $\A$-module $\A/\mathrm{ker}\, \varepsilon$ possesses a resolution
$(K_\ast,k)$ of the form
$$
	0 \rightarrow \A \rightarrow^{k_3} \A^3 \rightarrow^{k_2} \A^3
	\rightarrow^{k_1} \A \rightarrow \A / \ker \,\varepsilon \rto 0.
$$
The augmentation map $K_0=\A \rightarrow k = \A / \ker \,\varepsilon$ is given by the
counit $\varepsilon$. 
The left $\A$-module morphisms $k_n : K_n \rightarrow K_{n-1}$, 
$n=1,2,3$, are given 
by
$$
k_1 : 
\left(\begin{array}{c} 1 \\0 \\ 0 \end{array} \right),
\left(\begin{array}{c} 0 \\1 \\ 0 \end{array} \right),
\left(\begin{array}{c} 0 \\0 \\ 1 \end{array} \right) 
\mapsto a-1,b,c
$$
$$
k_2 : 
\left(\begin{array}{c} 1 \\0 \\ 0 \end{array} \right),
\left(\begin{array}{c} 0 \\1 \\ 0 \end{array} \right),
\left(\begin{array}{c} 0 \\0 \\ 1 \end{array} \right) 
\mapsto
\left(\begin{array}{c}
	b\\1-q^{-1}a\\0\end{array}\right),
	\left(\begin{array}{c}
	c\\0\\1-q^{-1}a\end{array}\right),
	\left(\begin{array}{c} 0 \\ c\\ -b\end{array}\right),
$$
\begin{eqnarray}
 k_3 : 1 \mapsto 
	\left(\begin{array}{c} c \\ -b \\ q^{-2}a-1 \end{array}\right).\nonumber
\end{eqnarray}  
\end{prop}  
\begin{pf}
It follows from Proposition~\ref{serrek}
(or directly) that this is a complex.
Let $(x,y,z)^t \in \ker (k_1)$, i.e.
$x(a-1)+yb+zc=0$. 
Then $\pi_{bc}(x(a-1))=0$. Using the $\bZ$-grading we have 
$\pi_{bc}(x)=0$, so $x=x'b+x''c$. Subtracting 
$k_2(x',x'' ,0)^t$ from $(x,y,z)^t$
we get a new element of $\ker (k_1)$ with $x=0$. Hence 
$\pi_c(y)=\pi_b(z)=0$, so $y=y'c$, $z=z'b$, $z'=-y'$ and this
element is a multiple of $k_2(0,0,1)^t$. 
In a similar manner  
$$
	x \left(\begin{array}{c}
	b\\1-q^{-1}a\\0\end{array}\right)+
	y \left(\begin{array}{c}
	c\\0\\1-q^{-1}a\end{array}\right)+
	z   
	\left(\begin{array}{c} 0 \\ c\\ -b\end{array}\right) = 0
$$
implies $x=x'c$, $y=-x'b$, $x'(q^{-2}a-1)=z$ for some $x' \in \A$.
\end{pf}
\begin{cor}\label{coro}
If $n>3$, then $HH_n^\s(\A)=0$ for all automorphisms $\sigma$. 
\end{cor}

The morphism between the resulting short complex 
$\sA' \otimes_\A K_\ast$ and the
standard complex for $\mathrm{Tor}_n^\A(\sA', k)$
yielding an isomorphism in homology is given explicitly by:\\
\begin{enumerate}
\item{
The map $\varphi_0 : \A \rightarrow C_0=\A$ is the identity.\\
}
\item{
The map $\varphi_1 : \A^3 \rightarrow C_1=\A^{\otimes 2}$ is given by
\begin{equation}\label{phi_1}
	\left(\begin{array}{c} x \\ y \\ z \end{array}\right) \mapsto
	(x,a-1)+(y,b) + (z,c).   
\end{equation}
Since $d(x,1,1) = (x,1)$ for any $x$, we have
$[(x, a-1)] = [(x,a)]$ in $\mathrm{Tor}_1^\A(\sA',k)$.\\
}
\item{
The map $\varphi_2 : \A^3 \rightarrow C_2=\A^{\otimes 3}$
is given by
$$
	\left(\begin{array}{c} x \\ y \\ z \end{array}\right) \mapsto
	(x,b,a-1) - (x, q^{-1}a-1, b) +$$
\begin{equation}\label{phi_2}
+ (y,c,a-1) - (y, q^{-1}a-1, c) + (z,c,b) - (z,b,c)
\end{equation}\\
}
\item{
Finally, the map $\varphi_3 : \A \rightarrow C_3=\A^{\otimes 4}$ in the
complex for the $ \mathrm{Tor} $-groups is given by
$x \mapsto x \otimes v$,  
where
$$
	v = -(q^{-2}a-1, b, c) + (q^{-2}a-1, c , b) - (c, q^{-1}a-1, b) +
$$
\begin{equation}\label{phi_3}
	+(c, b, a-1) -
	(b,c,a-1) +
	(b, q^{-1}a-1, c).
\end{equation}
}
\end{enumerate}
One sees by direct computation that this is a morphism of complexes,
and by the comparison theorem (see \cite{weibel}, Theorem~{2.2.6})
this is a quasi-isomorphism.

\subsection{\sc Computation of $HH_n^\s(\A)$, $n \leq 3$}

All automorphisms arising from finite-dimensional calculi are of the
form $\sigma=\sigma_{\lambda, \mu}$, and from now on we will only
consider automorphisms of this type. In fact, they are of the form 
$ \sigma(x)=$ $\sigma_{mod}^{-1} (f \ast x)$, where $f$ is a functional
in the dual Hopf algebra $\A^\circ$ acting on $x$ by 
$f \ast x=\sum f(x_{(2)})x_{(1)}$ (see Theorems~{4.1},~{4.3} 
and~{4.8} in \cite{kmt}). It is clear that such automorphisms
do not exchange $b$ and $c$.
By Corollary \ref{crucial_isomorphism} we have 
$HH_n^\s (\A) \cong H_n (\A, \sA)$, and the homologies 
$H_* (\A,\sA)$ can be calculated via our noncommutative Koszul resolution.

So let $\lambda,\mu \in k^\times$ and 
$\sigma=\sigma_{\lambda, \mu}$.
We apply $\sA'  \otimes_\A \cdot$ to our resolution and obtain the
complex $(F_\ast,f)$ of vector spaces
\begin{equation}
\label{complex}
	0 \rightarrow \A \rightarrow^{f_3} \A^3 \rightarrow^{f_2} 
	\A^3 \rightarrow^{f_1} \A \rightarrow 0,
\end{equation}
with morphisms $f_n$ given by
\begin{eqnarray}
&&
 f_1 \left(\begin{array}{c} x_1 \\ y_1 \\ z_1 \end{array}\right) 
	= x_1 \blacktriangleleft a- x_1+
	y_1 \blacktriangleleft b+ z_1 \blacktriangleleft c,\nonumber\\
&&
 f_2 \left(\begin{array}{c} x_2 \\ y_2 \\ z_2 \end{array}\right)  
	= \left(\begin{array}{c}
	x_2 \blacktriangleleft b + y_2 \blacktriangleleft c\\
	x_2 - q^{-1}x_2 \blacktriangleleft a + z_2 \blacktriangleleft c\\
	y_2 - q^{-1}y_2 \blacktriangleleft a - z_2 \blacktriangleleft b 
	\end{array}\right),\nonumber\\
&&
\label{defn_f_3}
\label{defn_f_2}
\label{defn_f_1}
 f_3 (x_3 ) = 
	\left(\begin{array}{c} x_3 \blacktriangleleft c 
	\\ -x_3 \blacktriangleleft b \\ q^{-2}x_3
	\blacktriangleleft a-x_3 \end{array}\right).
\end{eqnarray}  
Writing $\varepsilon_{i,j,k}:=q^{i+j+k+2} \lambda \mu^{-1}$
 we have 
\begin{eqnarray}\label{sindwiederda}
\label{right_action_of_a_b_c}
	\lambda q^{j+k} e_{i,j,k} \blacktriangleleft a 
&=& e_{i,j,k} + q^{-1-i-|i|}(1-\varepsilon_{|i|,j,k}) e_{i,j+1,k+1},\nonumber\\
	\lambda^{-1} e_{i,j,k} \blacktriangleleft b 
&=&	(1-\varepsilon_{i,j,k}^{-1}) 
	e_{i+1,j+1,k} \nonumber\\
&&+\left\{
	\begin{array}{ll}
	0\quad & : i \ge 0\\
	q^{-2i-1}(1-\varepsilon_{-i,j,k}^{-1})
	e_{i+1,j+2,k+1} \quad & : i<0
	\end{array}\right.,\nonumber\\
	\lambda e_{i,j,k} \blacktriangleleft c 
&=& (1-\varepsilon_{-i,j,k})e_{i-1,j,k+1} \nonumber\\
&&+\left\{
	\begin{array}{ll}
	q^{-2i+1}(1-\varepsilon_{i,j,k})
	e_{i-1,j+1,k+2} \quad & : i>0\\
	0 \quad & : i \le 0
	\end{array}\right.
\end{eqnarray}
For $q=\lambda=\mu=1$ we have $f_n=0$ and we recover the 
Hochschild-Kostant-Rosenberg theorem for $SL(2,k)$.
 The cotangent bundle of an algebraic group is trivial, so in the classical case 
$HH_n(\A)=\A \otimes \Lambda^n k^3$. It is clear, however, that for 
$q \neq 1$ there is no twisting automorphism for which this happens.

The calculations lead to five distinct cases:

\begin{enumerate}
\item{
\label{case_1} %Old case 1, 2, 4, 5
$\mu=1$, $\l \notin \{ q^{-(N+2)} \}_{ N \geq 0}$, and $\mu \neq 1$, $\l=1$.
}

\item{
\label{case_2} %Old case 3.
$\mu =1$, $\l = q^{-(N+2)}$, $N \geq 0$.
}

\item{
\label{case_3} %Old case 6
$\mu = q^{M+1}$, $\l=q^{-(N+1)}$, $M$, $N \geq 0$.
}

\item{
\label{case_4} %Old case 8
$\mu = q^{-(M+1)}$, $\l = q^{-(N+1)}$, $M$, $N \geq 0$.
}

\item{
\label{case_5} %Old cases 7 and 9
$\mu = q^{\pm(M+1)}$, $M \geq 0$, $\l \notin q^{-\bN}$, and $\mu \notin q^\bZ$, $\l \neq 1$.
}
\end{enumerate}

The computation of $HH_0^\s(\A)$ and $HH_1^\s(\A)$ 
is done most easily ``by hand" using the original Hochschild complex, but for 
$HH_2^\s(\A)$ and $HH_3^\s(\A)$ the calculations are done via the
Koszul resolution.

\subsection{$HH_0^\s (\A)$} 
We calculate from first principles the twisted Hochschild homology
$HH_0^\s (\A)$ for all automorphism $\s = \s_{\l,\mu}$. 
We start with the observation that: 
$$
	b(a_1,a_2a_3)=
	b(a_1a_2,a_3)+
	b(\sigma(a_3)a_1,a_2) \quad \forall \; a_1,\, a_2, \, a_3 \in\A
$$
So for any $a_1$, $a_2 \in\A$, there exist 
$x_a$, $x_b$, $x_c$, $x_c \in \A$ such that 
$$
	b( a_1 , a_2 ) = b[ (x_a, a) + (x_b ,b) + (x_c, c) + (x_d ,d)]
$$
Hence the image of the twisted Hochschild boundary is spanned by
\begin{eqnarray}
	A_{i,j,k}
&:=& e_{i,j,k}a - \lambda a e_{i,j,k} \nonumber\\ 
&=& (q^{-(j+k)}- \lambda) e_{i+1,j,k} \nonumber\\ 
&& + \left\{\begin{array}{ll}
	 0 & : i \ge 0\\
	(q^{-(j+k+1)}- \lambda q^{-2i-1}) e_{i+1,j+1,k+1}  & : i<0
	\end{array}\right. \nonumber\\ 
	B_{i,j,k}
&:=& e_{i,j,k} b - \mu b e_{i,j,k} = (1 - \mu q^{-i}) e_{i,j+1,k},\nonumber\\
	C_{i,j,k}
&:=&e_{i,j,k} c - \mu^{-1} c e_{i,j,k} =	(1-\mu^{-1} q^{-i}) e_{i,j,k+1},\nonumber\\
	D_{i,j,k}
&:=& e_{i,j,k} d - \lambda^{-1} d e_{i,j,k} \nonumber\\  
&=& (q^{j+k}- \lambda^{-1}) e_{i-1,j,k} \nonumber\\
&& +	 \left\{\begin{array}{ll}
	 0 \quad & : i \le 0\\
	(q^{j+k+1}- \lambda^{-1}q^{-2i+1}) e_{i-1,j+1,k+1} \quad & :  i>0
	\end{array}\right..\nonumber
\end{eqnarray}
For given $(i,j,k)$, the elements $B_{i,j,k}$ and $C_{i,j,k}$ both vanish
 if and only if $i=0$ and $\mu=1$. Therefore, for all $\lambda$, $\mu$,
$\mathrm{im}\, b$ contains the basis elements
\begin{equation}\label{diesinddrin}
	e_{i,j,k},\quad i \neq 0,\; j,k>0.
\end{equation} 
Omitting the span of these terms from the above list 
of generators we see that $\mathrm{im}\, b$ is spanned by 
(\ref{diesinddrin}) together with
\begin{eqnarray}
&& A_{-1,j,k} = (q^{-(j+k)}- \lambda) e_{0,j,k} +
	(q^{-(j+k+1)}- \lambda q) e_{0,j+1,k+1} 
	\nonumber\\ 
&& \tilde A_{i,j,k} = (q^{-(j+k)}- \lambda) e_{i+1,j,k},\quad
	i \neq -1 \nonumber\\
&& B_{i,j,0} = (1 - \mu q^{-i}) e_{i,j+1,0},\quad
	B_{0,j,k} = (1 - \mu ) e_{0,j+1,k},\nonumber\\
&& C_{i,0,k} =(1-\mu^{-1} q^{-i}) e_{i,0,k+1},\quad
	C_{0,j,k} = (1-\mu^{-1}) e_{0,j,k+1},\nonumber\\
&& D_{1,j,k} = (q^{j+k}- \lambda^{-1}) e_{0,j,k} +
	(q^{j+k+1}- \lambda^{-1}q^{-1}) e_{0,j+1,k+1} \nonumber\\
&& \tilde D_{i,j,k} = (q^{j+k}- \lambda^{-1}) e_{i-1,j,k},\quad
	i \neq 1.\nonumber
\end{eqnarray}
Since $\tilde D_{i+2,j,k}$ is proportional to $\tilde A_{i,j,k} $ and
both vanish if and only if $\lambda = q^{-(j+k)}$, we can omit 
$\tilde D_{i,j,k}$ from this list. We also have
$$
	A_{-1,j,k} = - \lambda q^{-(j+k)} D_{1,j,k},
$$
so the $D_{1,j,k}$ can be omitted as well. Finally, 
$C_{0,j,k}$ is for $j>0$ a nonzero multiple of 
$B_{0,j-1,k+1}$ and can be omitted. 
Thus the degree $0$ part (with respect to the 
$\bZ $-grading) of $\mathrm{im}\, b$
is spanned by 
$$(1- \lambda q^{j+k}) e_{0,j,k} + (q^{-1}- \lambda q^{j+k+1}) e_{0,j+1,k+1},$$ 
$$(1 - \mu ) e_{0,r,s},\quad r+s>0$$
and the nonzero degrees by (\ref{diesinddrin}) together with
$$
	(1- \lambda q^{j+k}) e_{i,j,k},\quad
	(1 - \mu q^{-i}) e_{i,j+1,0},\quad
	(1-\mu^{-1} q^{-i}) e_{i,0,k+1}.
$$
where $i \neq 0$, and $j$, $k \geq 0$.
Dually, 
$$HH^0_\s (\A) = \{ \; \linear \; h : \A \rto k \; : \; h(a_1 a_2 ) = h( \s(a_2 ) a_1 ) \; \}$$ 
For 
$\l \notin q^{-\bN}$ we have  $h( e_{i,j,k}) =0$ for $i \neq 0$, and
\begin{equation}\label{cohom_lambda_notin_q^-N}
 	h ( b^j c^k ) =  \left\{
	\begin{array}{ll}
	(-q)^{-k} {\frac{f(j-k)}{ f(j+k)}} h( b^{j-k} ) :& j \geq k\\
	{} & {}\\
	(-q)^{-j} {\frac{f(k-j)}{f(j+k)}}  h( c^{k-j}) :& j \leq k
	\end{array}\right. 
\end{equation}
where $f(n) = \l - q^{-n}$.

We now present the generating twisted 0-cycles, together with dual twisted 0-cocycles.
Our calculations now break down into five cases:\\

{\bf Case 1:} %Old cases 1,2,4,5
$\mu=1$, $\l \notin \{ q^{-(N+2)} \}_{N \geq 0}$ and $\mu \neq 1$, $\l=1$. Then
\begin{equation}
\label{HH_0_case_1}
HH_0^\s (\A) = k[1] \oplus \bigoplus_{x \in \{ a,b,c,d \}, \; \s(x) = x} ( \; \sum_{r \geq 0}^\oplus \; k [ x^{r+1} ] \; )
\end{equation}
For $\mu = 1 = \l$ (i.e. $\s = \id$) this agrees with \cite{mnw}. The dual $\s$-twisted 0-cocycles are defined on basis elements $x = e_{i,j,k}$ with $\s(x) =x$ as follows:
\begin{equation}
\label{h_1}
h_{[ 1 ]} (x) = \left\{
	\begin{array}{ll}
	1 & : x=1\\
	(-q)^{s+1} {\frac{f(0)}{f(2j+2)}} & : x = (bc)^{j+1}\\
	0 &: \otherwise
	\end{array}\right. 
\end{equation}
(if $\l =1$, obviously $f(0)=0$).
For $y = [ a^{r+1}]$, $[ d^{r+1}]$ define 
\begin{equation}\label{h_1_a^r+1_d^r+1}
	h_{[y]} (x) = \left\{
	\begin{array}{ll}
	1 & : x=y\\
	0 &: \otherwise
	\end{array}\right. 
\end{equation}
For $y = [ b^{s+1}]$, $[ c^{t+1}]$ define
\begin{equation}\label{h_b^s+1}
	h_{[ b^{s+1} ]} (x) = \left\{
	\begin{array}{ll}
	(-q)^k {\frac{f(s+1)}{ f(s+1+2k)}} & : x= b^{s+1} (bc)^k\\	
	0 &: \otherwise
	\end{array}\right. 
\end{equation}
\begin{equation}\label{h_c^t+1}
	h_{[ c^{t+1} ]} (x) = \left\{
	\begin{array}{ll}
	(-q)^j  {\frac{f(t+1)}{f(t+1+2j)}}   & : x= (bc)^j c^{t+1}\\
	0 &: \otherwise
	\end{array}\right. 
\end{equation}
These all satisfy  (\ref{cohom_lambda_notin_q^-N}).
For any $[x]$, $[y]$ in (\ref{HH_0_case_1}), we have
$h_{[y]} ( x ) = \d_{[x],[y]}$,  
so the 0-cycles given in (\ref{HH_0_case_1}) are linearly independent, hence a basis.\\

{\bf Case 2:} %This is old case 3.
$\mu =1$, $\l = q^{-(N+2)}$, $N \geq 0$. We  have 
\begin{equation}
\label{HH_0_case_2}
HH_0^\s (\A) = 
( \Sigma_{s \in S}^\oplus  \; k[ b^s ] ) \oplus 
( \Sigma_{t \in S}^\oplus  \; k[ c^t ] ) \oplus 
( \Sigma_{0 \leq i \leq N+2}^\oplus  \; k[ b^i c^{N+2-i} ] )
\end{equation}
where $S = \{ \integers \geq N+3 \} \cup \{ N+1, N-1, N-3, \ldots \geq 0 \}$, with the convention that if $0 \in S$, we include only one copy of $k[1]$.
Dual 0-cocycles are  $h_{[y]}$, defined for $[y] = [b^i c^{N+2-i}]$ on the basis $e_{i,j,k}$ by 
\begin{equation}
\label{h_b^i_c^N+2-i}
h_{[y]} (x) = \left\{
	\begin{array}{ll}
	1 & : x=y\\
	0 &: \otherwise
	\end{array}\right. 
\end{equation}
and  for $[y] = [ b^s]$, $[ c^t]$, $s$, $t \in S$ by
$$h_{[ b^s ]} (x) = \left\{
	\begin{array}{ll}
	(-q)^k {\frac{f(s)}{ f(s+2k)}} & : x= b^s (bc)^k\\	
	0 &: \otherwise
	\end{array}\right. $$
\begin{equation}
\label{h_b^s}
\label{h_c^t}
h_{[ c^t ]} (x) = \left\{
	\begin{array}{ll}
	(-q)^j  {\frac{f(t)}{f(t+2j)}}   & : x= (bc)^j c^t\\
	0 &: \otherwise
	\end{array}\right. 
\end{equation}
So for each pair $[x]$, $[y]$ appearing in (\ref{HH_0_case_2}), we have $h_{[y]} ( x ) = \d_{ [x], [y]}$.\\

{\bf Case 3:} %Old case 6
$\mu = q^{M+1}$, $\l = q^{-(N+1)}$, $M$, $N \geq 0$.
 \begin{equation}
\label{HH_0_oldcase_6}
HH_0^\s (\A) \cong k^2 = k[ d^{M+1} c^{N+1} ]  \oplus k [ a^{M+1} b^{N+1} ]
\end{equation}
Also $HH^0_\s (\A) \cong k^2$, with basis the twisted 0-cocycles $h_{[y]}$, 
$[y] = [ d^{M+1} c^{N+1} ]$, $[ a^{M+1} b^{N+1} ]$, defined on elements $x=e_{i,j,k}$, with $\s(x)=x$  by 
\begin{equation}\label{case_6_HH^0}
	h_{[y]} (x) = \left\{
	\begin{array}{ll}
	1 & : x=y\\
	0 &: \otherwise
	\end{array}\right. 
\end{equation}

{\bf Case 4:} %Old case 8
$\mu = q^{-(M+1)}$, $\l = q^{-(N+1)}$, $M$, $N \geq 0$.
We have 
\begin{equation}
HH_0^\s (\A) \cong k^2 = 
k[ d^{M+1} b^{N+1} ] \oplus k [ a^{M+1} c^{N+1} ]
\end{equation}
with $HH^0_\s (\A) \cong k^2$ with basis  $h_{[y]}$, $[y] = [ d^{M+1} b^{N+1} ]$, $[ a^{M+1} c^{N+1} ]$, defined as in (\ref{case_6_HH^0}).\\

{\bf Case 5:} %This is old cases 7 and 9.
$\mu = q^{\pm( M+1)}$, $M \geq 0$, $\l \notin q^{- \bN}$, and $\mu \notin q^{\bZ}$, $\l \neq 1$. Then
$$HH_0^\s (\A) =0 = HH^0_\s (\A)$$

\subsection{Twisted cocycles defined by derivations}
Before proceeding with $HH_1^\s(\A)$, we  
present a general construction 
of $\s$-twisted Hochschild $n$-cocycles.
It is essentially a variant of 
the characteristic map of \cite{como,cup}.

Let $\A$ be a $k$-algebra and $\s$ be an automorphism. 
A $\s$-derivation of $\A$ is a $k$-linear map 
$\der : \A \rto \A$, such that 
$\der( a_0 a_1 ) = a_0 \der (a_1 ) + \der( a_0 ) \s(a_1)$.
The following is straightforward:

\begin{prop}
\label{hoch_n_cocycle}
If $h$ is a $\s_1$-twisted 0-cocycle, $\der_1$, ..., $\der_{n-1}$ derivations of 
$\A$, and $\der_n$ is a $\s_2$-derivation, then  
\begin{equation}
\label{defn_phi_n}
	\phi_n ( a_0, a_1, \ldots , a_n ) = 
	h( a_0 \der_1 (a_1) \ldots \der_n (a_n) )
\end{equation}
is a $\s_1 \circ \s_2$-twisted Hochschild n-cocycle 
\end{prop}

In general there is no reason for such a cocycle to represent a nontrivial element of Hochschild cohomology, nor for it also to be cyclic. However:

\begin{lemma}
\label{cyclic_1_cocycle}
 Suppose $\A$ is a unital algebra, $h$ a $\s_1$-twisted $0$-cocycle and 
$\der$ a $\s_2$-derivation of $\A$.
Defining $\phi_1$ by $\phi_1 (x,y) = h( x \der (y))$, 
then $\phi_1$ is a $\s_1 \circ \s_2$-twisted cyclic cocycle 
if and only if $h ( \der(a) ) = 0$ for all $a \in \A$.
\end{lemma}

For $\A = \A( \slq2)$ there are obvious derivations $\der_a$, $\der_b$ defined by
$$\der_a (a) = a, \quad
\der_a (b) = 0, \quad
\der_a (c) = 0, \quad
\der_a (d) = -d$$
\begin{equation}
\label{der_a}
\label{der_b}
\der_b (a) = 0, \quad
\der_b (b) = b, \quad
\der_b (c) = -c, \quad
\der_b (d) = 0
\end{equation}
and extended via the Leibniz rule. 
%Also, corresponding to 
For any $x \in \A$ define an inner derivation $\der'_x$  by 
$\der'_x (y) = [x,y] = xy-yx$. The following is straightforward:

\begin{prop} 
The vector space of all derivations of $\A ( \slq2)$ is spanned by $\der_a$, $\der_b$ together with the inner derivations.
\end{prop}

In the sequel we will use the derivation $\der_0 = \der_a + \der_b$, which satisfies
\begin{equation}
\label{der}
\label{der_0}
\der_0 (a) = a, \quad
\der_0 (b) = b, \quad
\der_0 (c) = -c, \quad
\der_0 (d) = -d
\end{equation}
and also the $\s_{\l,1}$-derivation defined by
\begin{equation}
\label{sigma_der}
\der(a) = a, \quad \der(b)= 0 = \der(c), \quad \der(d) = - \l^{-1} d.
\end{equation}
 
\subsection{$HH_1^\s (\A)$}

The second twisted Hochschild boundary is given by
$$
	b(a_1,a_2,a_3)=
	(a_1a_2,a_3)-(a_1,a_2a_3)+
	(\sigma(a_3)a_1,a_2).
$$
In particular (take $a_2=a_3=1$) the image of $b$ contains all elementary tensors 
of the form $(a_1,1)$, and the residue classes of 
$$
	(e_{i,j,k},a),\quad
	(e_{i,j,k},b),\quad
	(e_{i,j,k},c),\quad
	(e_{i,j,k},d)
$$ 
generate $\A \otimes \A / \mathrm{im}\, b$. Now,
$HH_1^\s (\A)$ is the kernel of the map 
$\A \otimes \A / \mathrm{im}\, b \rightarrow \A$
 induced by the first
twisted Hochschild boundary. This sends the classes of the above 
elements to $A_{i,j,k},B_{i,j,k},C_{i,j,k},D_{i,j,k}$ from the
previous section. It is straightforward to check for triviality and
linear dependence. 
%In doing so one arrives at the following presentations of $HH_1^\s (\A)$:

We now present generators of $HH_1^\s (\A)$ 
and dual twisted 1-cocycles. 
From Proposition \ref{hoch_n_cocycle}, 
for any $\s$-twisted 0-cocycle $h$ and 
derivation $\der$,  defining $\phi_1(x,y) = h( x \der(y) )$ gives a $\s$-twisted Hochschild 1-cocycle. 
For $\A(\slq2)$ all automorphisms $\s_{\l,\mu}$ commute with the
derivations $\der_a$, $\der_b$ (\ref{der_a}), and the 
$e_{i,j,k}$ are eigenvectors for these derivations. 
By Lemma \ref{cyclic_1_cocycle},  $\phi_1$ is
cyclic if and only if $h \circ \der =0$.  
We take $\der_0 = \der_a + \der_b$ defined by (\ref{der}).\\

{\bf Case 1:} $\mu=1$, $\l \notin \{ q^{-(N+2)} \}_{N \geq 0}$ and $\mu \neq 1$, $\l=1$. Then
\begin{equation}
\label{HH_1_case_1}
HH_1^\s (\A) = k[\omega_1] \; \oplus \bigoplus_{x \in \{ a,b,c,d \}, \; \s(x) = x} ( \; \sum_{r \geq 0}^\oplus k [ (x^r ,x)] ) 
\end{equation}
where $\om_1 = (\mu^{-1} -1) (d,a) + (q - q^{-1})(b,c)$. 
We note that, for all $\mu$ and for $\l \neq q^{-2}$, we have
$[(c,b)] = - \mu [(b,c)]$, $[(d,a)] = - \l [(a,d)]$.
For $\mu = 1 =\l$ this is in agreement with \cite{mnw}, apart from the sign change in $\om_1$.
Now recall the
0-cocycles $h_{[x]}$ defined in (\ref{h_1})-(\ref{h_c^t+1}).
%, (\ref{h_1_a^r+1_d^r+1}), (\ref{h_b^s+1}), (\ref{h_c^t+1}). 
Given such an $x$, define a Hochschild
1-cocycle $\phi_{[x]}$ by 
\begin{equation}\label{phi_[x]}
	\phi_{[x]} ( y,z) = h_{[x]} ( y \der_0 (z) )
\end{equation}
Then the Hochschild 1-cocycles dual to the generators of $HH_1^\s (\A)$ are:
$$
		\phi_{[ x^{r+1} ]} \leftrightarrow ( x^r ,x), \quad 
		x \in \{ a,b,c,d \}, \quad \s(x) = x.
$$
Dual to $ \omega_1 $ we have 
\begin{equation}\label{phi_omega_1}
		\phi_{\omega_1} (x,y) =
		h_{[1]} (x \der_0(y))
\end{equation}
with $h_{[1]}$ defined in (\ref{h_1}) and
$\der_0$ defined in (\ref{der_0}). Then for $\l=1$ we have
$ \phi_{\omega_1}(d,a)=1$, and for $\l \neq q^{-2}$ we have 
 $ \phi_{\omega_1}(b,c)=\frac{q(1-\l)}{\l q^2-1}$.
 Since $h_{[1]} \circ \der_0 =0$,  by Lemma \ref{cyclic_1_cocycle}  $\phi_{\om_1}$ is in fact a $\s$-twisted cyclic 1-cocycle.\\

{\bf Case 2.} %Old case 3
$\mu =1$, $\l =q^{-(N+2)}$, $N \geq 0$. Then 
$$HH_1^\s (\A) = 
( \Sigma_{s \in S'}^\oplus \; k[( b^s ,b)] ) \oplus
( \Sigma_{t \in S'}^\oplus \; k[( c^t ,c)] ) \oplus$$
\begin{equation}
\label{HH_1_case_3a}
( \Sigma_{0 \leq i \leq N}^\oplus \; k[( b^i c^{N+1-i} ,b )] ) \oplus
( \Sigma_{0 \leq i \leq  N}^\oplus \; k[( b^{i+1} c^{N-i} ,c )] ) \oplus
k [ \omega_1 ]
\end{equation}
Here $S' = \{ \integers \geq N \} \cup \{ N-2, N-4, \ldots \geq 0 \}$, and for $N \geq 1$ we have 
 $[ \omega_1] = [(b,c)] = - [(c,b)]$ for $N$ odd, $[\omega_1]=0$ for $N$ even. For $\l = q^{-2}$, $[(b,c)]$ and $[(c,b)]$ are linearly independent.

Recall the 0-cocycles $h_{[x]}$, with
$[x] =[b^s]$, $[c^t]$, defined in 
(\ref{h_b^s+1}), (\ref{h_c^t+1}). 
The dual Hochschild 1-cocycles (\ref{phi_[x]}) are
$$\phi_{[b^{s+1}]} \leftrightarrow (b^s ,b), \quad
\phi_{[c^{t+1}]} \leftrightarrow (c^t ,c) $$
together with the twisted cyclic 1-cocycle $\phi_{\om_1}$ (\ref{phi_omega_1}) dual to $\om_1$.
 To define twisted 1-cocycles dual to 
$[( b^i c^{N+1-i} ,b )]$, $[( b^{i+1} c^{N-i} ,c )]$ we need to work a little harder.
It is straightforward to show that any $\s$-twisted Hochschild 1-cocycle is uniquely defined by its values $\phi(y,t)$, for $t = a,b,c,d$ and basis elements $y = e_{i,j,k}$.

\begin{lemma}
For $\l = q^{-(N+2)}$, $\mu=1$ defining $\phi_1$, $\phi_2$ 
on basis elements $y=e_{i,j,k}$ by
%$$
%	\phi_1 ( b^i c^{N+1-i} ,b) = 
%	\b_{N,i}, \quad \phi_1 ( d b^i c^{N-i} , a) = 
%	q^{-(N+1)} \b_{N,i},
%$$
\begin{eqnarray}
&& \phi_1 ( b^i c^{N+1-i} ,b) = \b_{N,i}, \quad \phi_1 ( d b^i c^{N-i} , a) = q^{-(N+1)} \b_{N,i},\nonumber\\
&& \phi_1 (y,b) = 0 = \phi_1 (y,a) \quad \otherwise, \quad  \phi_1 (y,c) = 0 = \phi_1 (y,d) \quad \forall \; y\nonumber\\
&& \phi_2 ( a b^i c^{N-i} , d) = q^{N+1} \g_{N,i}, \quad \phi_2 ( b^{i+1} c^{N-i} , c) = \g_{N,i} \nonumber\\
&& \phi_2 (y,d) = 0 = \phi_2 (y,c) \quad \otherwise, \quad  \phi_2 (y,a) = 0 = \phi_1 (y,b) \quad \forall \; y\nonumber
\end{eqnarray}
for arbitrary $\b_{N,i}$, $\g_{N,i}$, $0 \leq i \leq N$,
gives well-defined $\s$-twisted Hochschild 1-cocycles.
\end{lemma}

Setting each $\b_{N,i}$, $\g_{N,i}$ to 1 in turn, and all  others to zero, we see that the twisted Hochschild 1-cycles $[( b^i c^{N+1-i} ,b )]$, $[( b^{i+1} c^{N-i} ,c )]$  are nontrivial and linearly independent.
This extends to $N+1$ linearly independent twisted cyclic 1-cocycles. Define $\phi = \phi_1 + \phi_2$. Cyclicity requires
%$$\phi( b^i c^{N+1-i} ,b) = - \phi( b, b^i c^{N+1-i}) = -i \phi( b^i c^{N+1-i} ,b) - (N+1-i) \phi ( b^{i+1} c^{N-i} ,c)$$
$$\phi( xc ,b) = - \phi( b, xc) = -i \phi( xc,b) - (N+1-i) \phi ( bx ,c)$$
where $x = b^i c^{N-i}$.
Hence for $\phi$ to be cyclic, we need
$$(i+1) \b_{N,i} = -(N+1-i) \g_{N,i} \quad 0 \leq i \leq N$$
i.e. $\g_{N,i} = {\frac{-(i+1)}{N+1-i}} \b_{N,i}$ for $0 \leq i \leq N$. Then:

\begin{lemma}
For each  $0 \leq i \leq N$, defining
 $\phi_{N,i} = \phi_1 + \phi_2$ with 
 $$\b_{N,i} = N+1-i, \quad \g_{N,i} = -(i+1)$$
 gives a well-defined twisted cyclic 1-cocycle satisfying
 \begin{eqnarray}
&& \phi_{N,i} ( xc ,b) = N+1-i,\quad \phi_{N,i} ( bx ,c) = -(i+1)\nonumber\\
&& \phi_{N,i} (ax,d) = -q^{N+1} (i+1), \quad \phi_{N,i} (dx,a) = q^{-(N+1)} (N+1-i)\nonumber
\end{eqnarray}
for $x= b^i c^{N-i}$, and $\phi_{N,i} (y,t) = 0$ for all basis elements $y$ and $t=a,b,c,d$ otherwise.\\
\end{lemma}
%We will see in Section \ref{section_tc_slq2} that the  $\phi_{N,i}$  are a basis for  $HC^1_\s (\A) \cong k^{N+1}$.\\

{\bf Case 3.} %Old case 6
$\mu = q^{M+1}$, $\l=q^{-(N+1)}$, $M$, $N \geq 0$. Then
$HH_1^\s (\A)  \cong k^4$, with basis given by
 the Hochschild cycles
\begin{equation}
\label{case6_HH1}
( a^M b^{N+1} ,a), \quad
( a^{M+1} b^N ,b), \quad
( d^{M+1} c^N ,c), \quad
( d^M c^{N+1} ,d)
 \end{equation}
 The dual basis for $HH^1_\s (\A) \cong k^4$ is given by the 1-cocycles 
 $\phi_{[x],t}$, for $[x] = [ a^{M+1} b^{N+1} ]$, $[d^{M+1} c^{N+1}]$, and $t = a, b$ defined by
 \begin{equation}
 \label{phi_[x]_t}
 \phi_{[x],t} (y,z) = h_{[x]} ( y \der_t (z) )
 \end{equation}
 where the $h_{[x]}$ were defined in (\ref{case_6_HH^0})  and $\der_a$, $\der_b$ in (\ref{der_a}).
 We have
 \begin{eqnarray}
&& \phi_{ [ a^{M+1} b^{N+1} ], a} ( a^M b^{N+1} ,a) = q^{-(N+1)} \nonumber\\ 
&& \phi_{ [ a^{M+1} b^{N+1} ],b} ( a^{M+1} b^N ,b) = 1 \nonumber\\
&& \phi_{[d^{M+1} c^{N+1}] ,b}( d^{M+1} c^N ,c) = -1 \nonumber\\
&& \phi_{[d^{M+1} c^{N+1}],a} ( d^M c^{N+1} ,d) = -q^{N+1} \nonumber
\end{eqnarray}
with all other pairings being zero.\\

{\bf Case 4.} %Old case 8.
$\mu = q^{-(M+1)}$, $\l = q^{-(N+1)}$, $M$, $N \geq 0$.
$HH_1^\s (\A) \cong k^4$, with basis given by the 
 Hochschild cycles
\begin{equation}
 ( a^M c^{N+1} ,a ), \quad 
 ( d^{M+1} b^N ,b), \quad
( a^{M+1} c^N , c), \quad
( d^M b^{N+1} ,d)
\end{equation}
Analogously to Case 3, the dual basis  for $HH^1_\s (\A) \cong k^4$ is given by the 1-cocycles 
 $\phi_{[x],t}$, (\ref{phi_[x]_t}) for $[x] = [ a^{M+1} c^{N+1} ]$, $[d^{M+1} b^{N+1}]$, and $t = a$, $b$.\\

{\bf Case 5.} For $\mu = q^{\pm(M+1)}$, $M \geq 0$, $\l \notin q^{-\bN}$, and 
$\mu \notin q^\bZ$, $\l \neq 1$,
 $HH_1^\s (\A) =0$ and $HH^1_\s (\A) =0$.\\

\subsection{\sc $HH_2^\s (\A)$}

We now compute $HH_2^\s (\A)$. The first step is to describe the kernel of
$$
	\psi : \A^2 \rightarrow \A,\quad (x,y) \mapsto \psi_+(x)+\psi_-(y)
$$
If $ \mathrm{ker}\, \psi_\pm^c \subset \A$ are fixed complements to 
$ \mathrm{ker}\, \psi_\pm$ and 
$$
	\phi_\pm : 
	\mathrm{im}\, \psi_\pm \rightarrow 
	\mathrm{ker}\, \psi_\pm^c
$$
are the inverses of $ \psi_\pm |_{\mathrm{ker}\, \psi_\pm^c}$, then
the linear map
\begin{eqnarray}
&& \mathrm{ker}\, \psi_+ \oplus 
	(\mathrm{im}\, \psi_+ \cap \mathrm{im}\, \psi_-) \oplus
	\mathrm{ker}\, \psi_- \rightarrow  
	\mathrm{ker}\, \psi,\nonumber\\ 
&& (x,y,z) \mapsto 
	(x+\phi_+(y),z-\phi_-(y)) \nonumber
\end{eqnarray} 
is an isomorphism of vector spaces. Using (\ref{sindwiederda}) 
one determines a basis of 
$\mathrm{im}\, \psi_+ \cap \mathrm{im}\, \psi_-$ and obtains:
\begin{prop}\label{eswird}
The set
\begin{eqnarray}
&& \B_{\mathrm{ker}\,\psi} :=	\{(x_+,0),(0,x_-)\,|\,x_\pm \in \B_{\mathrm{ker}\,\psi_\pm}\}\nonumber\\
&& \cup \{(-\lambda^{-2}(
	\frac{1-\varepsilon_{-i,j,k}}{1-\varepsilon_{i-2,j,k}^{-1}}e_{i-2,j-1,k+1}
	-\varepsilon_{-i,j,k+1}e_{i-2,j,k+2}),e_{i,j,k}),\nonumber\\
&& \quad(e_{-i,k,j},-\lambda^2(
	\frac{1-\varepsilon_{-i,j,k}^{-1}}{1-\varepsilon_{i-2,j,k}}e_{-i+2,k+1,j-1}
	-\varepsilon_{-i,j,k+1}^{-1}e_{-i+2,k+2,j}))\,|\,\nonumber\\ 
&& \quad i \ge 2,j \ge 1,k \ge 0,\varepsilon_{i-2,j,k} \neq 1\} 
	\nonumber\\
&&\cup \{(\lambda^{-2}qe_{i-2,0,j+2},e_{i,0,j}),
	(e_{-i,j,0},\lambda^2q^{-1} e_{-i+2,j+2,0})\,|\, 
	i \ge 2,j \ge 0,\varepsilon_{-i,0,j}=1\} 
	\nonumber\\
&& \cup \{(e_{-1,j,k},\lambda^2 \varepsilon_{-1,j,k}^{-1} e_{1,j+1,k-1})\,|\,
	j \ge 0,k \ge 1\}.\nonumber
\end{eqnarray} 
is a vector space basis of $\mathrm{ker}\, \psi$.
\end{prop}

Now we can compute which of these remain nontrivial
and linearly independent 
modulo the 
image of the map $\varphi : x \mapsto (\psi_-(x),-\psi_+(x))$. 
\begin{prop}
The classes of 
\begin{eqnarray}
&& \{ (e_{i,j,0},0), \; (0,e_{-i,0,j}) \; | \; i \ge 0,\; j \ge 0,\; \varepsilon_{i,0,j}=1\; \}
	\nonumber\\
&& \cup \; \{ (\lambda^{-2}qe_{i-2,0,k+2}, \; e_{i,0,k}) \; | \; 
	i \ge 2,\; k \ge 0,\; \varepsilon_{-i,0,k}=1\; \} 
	\nonumber\\
&& \cup \;  
	\{ (e_{i,j,0}, \; \lambda^2 q^{-1} e_{i+2,j+2,0}) \; | \; 
	i \le -2, \; j \ge 0, \; \varepsilon_{i,j,0}=1 \; \} 
	\nonumber\\
&& \cup \; \{ (e_{-1jk},\lambda^2 e_{1,j+1,k-1}) \; | \;
	j \ge 0,\; k \ge 1,\varepsilon_{-1,j,k}=1 \; \}\nonumber
\end{eqnarray}
form a vector space 
basis of $\mathrm{ker}\, \psi / \mathrm{im}\, \varphi$. 
\end{prop}

Next we check for which linear combinations 
$(x,y) \neq (0,0)$ of these 
there exists $z$ with $(x,y,z) \in \mathrm{ker}\, f_2$, 
and determine those $z$ with 
$(0,0,z) \in \mathrm{ker}\, f_2$, giving a
generating set for $HH_2^\s(\A)$:
\begin{prop}
The classes of
\begin{eqnarray}
&& \{\; (0,0,e_{0,j,k}) \; | \; j, k \ge 0,\; \varepsilon_{0,j,k}=1\; \}\nonumber\\
&& \cup \; \{\; (e_{i,j,0},0,0),\; (0,e_{-i,0,j},0) \; | \; i \ge 0,\; j \ge 0, \; 
	\lambda=q^{-j-1},\; \mu=q^{i+1} \; \}
	\nonumber\\
&& \cup \{ (\lambda^{-2} q e_{i-2,0,k+2}, e_{i,0,k},
	\lambda^{-1} q^{-1} e_{i-1,0,k+1}) \; | \; 
	i \ge 2,k \ge 0,\lambda=q^{-k-1},\mu=q^{-i+1} \} 
	\nonumber\\
&& \cup \{ (e_{i,j,0}, \; \lambda^2 q^{-1} e_{i+2,j+2,0},
	\lambda q^{-1} e_{i+1,j+1,0}) \; | \;  
	i \le -2,j \ge 0,\lambda=q^{-j-1},\mu=q^{i+1} \} 
	\nonumber\\
&& \cup \{\; (e_{-1,j,k},\; \lambda^2 e_{1,j+1,k-1},
	\lambda q^{-1} e_{0,j+1,k}) \; | \;
	j \ge 0,\; k \ge 1,\; \lambda=q^{-j-k-1},\; \mu=1\; \}\nonumber
\end{eqnarray}
generate $HH_2^\s(\A)$.
\end{prop} 
The classes of the elements in the first line are trivial for 
$\mu \neq 1$, and for $ \mu=1$ they contain those from the last line.  
It follows directly from the definition of $f_3$ that the 
remaining classes are independent. Hence: 
\begin{equation}\label{HH_2}
	HH_2^\s (\A) \cong \left\{
	\begin{array}{ll}
	{k}^{N+1} \quad & : 
	\lambda=q^{-(N+2)}, N \ge 0,\; \mu=1,\\
	{k}^2 \quad & : 
	\lambda=q^{-(N+1)}, \;
	\mu=q^{\pm (M+1)}, M, N \geq 0,\\
	0 \quad & : \otherwise.
	\end{array}\right.  
\end{equation}

To calculate $HC_\ast^\s(\A)$ we need generators in the original 
Hochschild complex. In Case~2 we compute generators  from the above using
$\varphi_2$ (\ref{phi_2}) and  
$\xi$ (\ref{fetsimo}). In Cases 3 and 4 we use simpler
generators that are directly verified to be homologous to those obtained
from the above ones:\\

{\bf Case 2.} %Old case 3
$\mu=1$, $\l = q^{-(N+2)}$, $N \geq 0$.
Then $HH_2^\s (\A) \cong k^{N+1}$. Taking $x = b^i c^{N-i}$ ($0 \leq i \leq N$), a basis is given by
$$
	\omega_2 (N,i) = (bcx,a,d) - (bcx,d,a) -q (dbx,a,c) + q (bdx,c,a)
$$
\begin{equation}\label{HH_2_case_3}\label{omega_2_N_i}
	+ (dax,b,c) -(adx,c,b) - q^{-1} (cax,b,d) + q^{-1} (acx,d,b)
\end{equation}\\

{\bf Case 3.} %Old case 6
$\mu = q^{M+1}$, $\l=q^{-(N+1)}$, $M$, $N \geq 0$.
Then $HH_2^\s (\A) \cong k^2$, with basis given by the 
Hochschild cycles
\begin{eqnarray}
\label{case6_w2}\label{case6_w2p}
&\omega_2 = &(a^M b^N ,b,a) - q^{-1} ( a^M b^N ,a,b) \nonumber\\
&\omega_2^{'} = &( d^M c^N ,c,d) - q( d^M c^N ,d,c)
\end{eqnarray}

{\bf Case 4.} %Old case 8!
$\mu = q^{-(M+1)}$, $\l = q^{-(N+1)}$, $M$, $N \geq 0$.
$HH_2^\s (\A) \cong k^2$, with basis given by 
 the Hochschild cycles
 \begin{eqnarray}
&\omega_2 = &( a^M c^N ,c,a) - q^{-1} ( a^M c^N ,a,c) \nonumber\\
&\omega_2^{'} = &(d^M b^N ,b,d) - q( d^M b^N ,d,b)
\end{eqnarray}

Finally, $HH_2^\s (\A)=0$ for all other $\s = \s_{\l,\mu}$.

\subsection{$HH_3^\s (\A)$}

The third homology $HH_3^\s (\A)$ can be determined easily 
using the Koszul resolution. We abbreviate:
$$
	\psi_+: \A \rightarrow \A,\quad x \mapsto x \blacktriangleleft b,
	\quad 
	\psi_- : \A \rightarrow \A,\quad x \mapsto x \blacktriangleleft c.
$$
From (\ref{right_action_of_a_b_c}) we obtain in a straightforward way:
\begin{prop}\label{ersterkern}
The sets 
$$
	\B_{\mathrm{ker}\,\psi_\pm} := 
	\{ \; e_{i,j,k} \; | \; \pm i \ge 0,\; \varepsilon_{\pm i,j,k}=1 \; \}
$$
are vector space bases of $\mathrm{ker}\,\psi_\pm$. Hence
the sets
$$
	\B_{\mathrm{im}\,\psi_\pm} := 
	\{ \; E^\pm_{i,j,k} \; | \; e_{i,j,k} \notin \mathrm{ker}\, \psi_\pm \; \},\quad
	E^\pm_{i,j,k}:=\psi_\pm(e_{i,j,k})
$$
are vector space bases of $\mathrm{im}\,\psi_\pm$.
\end{prop}

If $x \in HH_3^\s (\A)=\mathrm{ker}\,f_3$, then 
$x \in \mathrm{ker}\, \psi_+ \cap \mathrm{ker}\, \psi_-$, so
by Proposition~\ref{ersterkern}, 
$\pi_{i,j,k}(x) \neq 0$ implies $i=0$. Insertion in
$q^{-2}x \blacktriangleleft a=x$ gives
$$
	x \in 
	\left\{
	\begin{array}{ll}
	\mathrm{span}\, \{b^i c^{N-i}\}_{i=0,\ldots,N} \quad & :
	\lambda=q^{-N-2},\mu=1,\\
	0 \quad & : \otherwise.
	\end{array}\right.  
$$
Conversely, all these monomials are elements of $\mathrm{ker}\, f_3$.
Hence:\\

{\bf Case 2:} %This is old case 3.
For $\mu =1$, $\l = q^{-(N+2)}$, $N \geq 0$ we  have 
$HH_3^\s(\A) \cong k^{N+1}$.\\

{\bf Cases 1, 3, 4, 5:} %This is old case 3.
$HH_3^\s(\A)=0$.\\

It is also straightforward that $HH_3^\s (\A)=0$
for all $\s= \tau_{\lambda,\mu}$. Therefore:

\begin{thm}
\label{HH_3_for_any_sigma}
For any automorphism $\s$, we have
$$
	HH_3^\s (\A)=\left\{
	\begin{array}{ll}
	{k}^{N+1} \quad & : \s = \s_{\l, \mu}, \;
	\lambda=q^{-(N+2)},N \geq 0,\; \mu=1,\\
	0 \quad & : \otherwise.
	\end{array}\right.  
$$
\end{thm}

Note that the $N=0$ case ($\l = q^{-2}$, $\mu=1$)  is precisely the modular automorphism (\ref{smod_slq2}).
For $\l = q^{-(N+2)}$ we translate 
the generators back to the original Hochschild complex
using the maps $ \varphi_3$ (\ref{phi_3}) and 
$ \xi $ (\ref{fetsimo}), giving
\begin{equation}
\label{omega_3_N_i}
	\omega_3 (N,i) = A(N,i) - B(N,i),\quad
	0 \leq i \leq N,
\end{equation}
\begin{eqnarray}
&A(N,i) = &dx \otimes (a \wedge b \wedge c) 
+ cx \otimes (b \wedge a \wedge d) \nonumber\\
&B(N,i) = &-qdbx \otimes (1 \wedge a \wedge c) - q^{-1} cax \otimes (1 \wedge b \wedge d) \nonumber\\
&&+ dax \otimes (1 \wedge b \wedge c) + bcx \otimes ( 1 \wedge a \wedge d) \nonumber\\
&& + (q - q^{-1}) bcx \otimes ( (c,b,1) - (1,c,b) + (c,1,b)) \nonumber
\end{eqnarray}
with $x= b^i c^{N-i}$, and the terms $``a_0 \wedge a_1 \wedge a_2"$ are given by:
\begin{eqnarray}
&a \wedge b \wedge c =& (a,b,c) - (a,c,b) + q(c,a,b) - q^2(c,b,a) +
q^2 (b,c,a) - q(b,a,c) \nonumber\\
&b \wedge a \wedge d =& (b,a,d) - (b,d,a) + q (d,b,a) - (d,a,b) + (a,d,b) - q^{-1}(a,b,d)\nonumber\\
&1 \wedge a \wedge c =& (1,a,c) - q(1,c,a) + q (c,1,a)
- q(c,a,1) + (a,c,1) - (a,1,c)\nonumber\\
&1 \wedge b \wedge d =& (1,b,d) - q (1,d,b) - (b,1,d) + (b,d,1) - q(d,b,1) + q (d,1,b)\nonumber\\
&1 \wedge b \wedge c =& (1,b,c) - (1,c,b) - (b,1,c) + (b,c,1) + (c,1,b) - (c,b,1)\nonumber\\
&1 \wedge a \wedge d =& (1,a,d) - (1,d,a) + (d,1,a) - (d,a,1) + (a,d,1) - (a,1,d)\nonumber
\end{eqnarray}
and throughout we denote $a_0 \otimes a_1 \otimes a_2$ by 
$(a_0, a_1,a_2)$.

In the normalized complex this becomes 
$ \omega_3(N,i)=A(N,i)$ since $B(N,i)$ is degenerate.

\section{Twisted cyclic homology of $\A( SL_q (2))$}
\label{section_tc_slq2}

We calculate the twisted cyclic homology of $\A (SL_q (2) )$ as the
total homology of Connes' mixed $(b,B)$-bicomplex
(\ref{mixed_b_B_bicomplex}) coming from the underlying cyclic object, as
in section \ref{section:cyclic_module}. Having found the twisted
Hochschild homology, we can now complete the spectral sequence
calculation. We remind the reader that throughout we are working with the
normalized mixed complex.

\subsection{Case 1}

\begin{prop}
In case 1, $\mu = 1$, $\l \notin \{ q^{-(N+2)} \}_{N \geq 0}$, and $\mu \neq 1$, $\l=1$,  %old cases 1,2,4,5! 
$HC_0^\s(\A)$ is infinite-dimensional, 
$HC_{2n+1}^\s(\A) = k[ \omega_1]$, and
$HC_{2n+2}^\s(\A) = k[1]$,
where  $[\omega_1]$ is the distinguished generator of
$HH_1^\s (\A)$. 
\end{prop}
\begin{pf}
By definition, $HC_0^\s(\A) = HH_0^\s (\A)$,  generated by $[1]$, together with
$[x^{r+1}]$ ($r \geq 0$), for those $x \in \{a,b,c,d\}$  with $\s(x) = x$,
while $HH_1^\s (\A)$ is generated by $[(x^r, x)]$, $(r \geq 0$) for the same set of $x$, together with the distinguished generator $[\omega_1 ] = (\mu^{-1} -1) [(d,a)] + (q - q^{-1}) [(b,c)]$.
We have
\begin{eqnarray} 
&&B_0 [1] = [(1,1)]= [ b (1,1,1)]=0,\nonumber\\
&&B_0 [x^{r+1}] = [(1,x^{r+1})] = (r+1) [(x^r ,x)]\nonumber
\end{eqnarray}
 by Lemma \ref{lemma:B_0_x^s_y^t}.
Hence $\ker(B_0 ) = k[1]$, and $HH_1^\s (\A) = \im(B_0 ) \oplus k [\omega_1]$. 
Further $HH_n^\s (\A) = 0$ for $n \geq 2$ in each case.
Hence the spectral sequence stabilizes at the second page:

	\[
     	\begin{CD}
 	@ VVV @ VVV @ VVV @ . @ . @ . @ .\\
	{0} @ <<< {k [\omega_1]} @ <<< {k[1]} @ . @ . @ . @ . @ .\\
 	@ VVV @ VVV @ . @ . @ . @ . @ .\\
{k[\omega_1]} @ <<< {k[1]} @ . @ . @ . @ . @ . \\
 	@ VVV @ . @ .  @ . @ . @ . @ .\\
	{HH_0^\s (\A)} @ . @ . @ . @ . @ . @ .\\
    	\end{CD}
	\]\\

with all further maps being zero. The result follows. \end{pf}

\subsection{Case 2} %Old case 3!

\begin{prop}
\label{HC_case_2}
In case 2, $\mu =1$, $\l = q^{-(N+2)}$, $N \geq 0$, we have 
$HC_0^\s(\A)$ infinite dimensional, while
$$HC_1^\s (\A) \cong   
	\left\{
	\begin{array}{ll}
	 k^{N+1} &: N \; \even\\
	 k^{N+2} &: N \; \odd
	\end{array}\right., \quad HC_2^\s (\A) \cong k^{N+2}, \quad N \; \odd$$
%{\bf Higher $HC_n^\s (\A)$ are still not found!}
\end{prop}
\begin{pf}
Recall from (\ref{HH_0_case_2}) that 
$$HH_0^\s (\A) = 
( \Sigma_{s \in S}^\oplus  \; k[ b^s ] ) \oplus 
( \Sigma_{t \in S}^\oplus  \; k[ c^t ] ) \oplus 
( \Sigma_{0 \leq i \leq N+2}^\oplus  \; k[ b^i c^{N+2-i} ] )$$
where $S = \{ \integers \geq N+3 \} \cup \{ N+1, N-1, N-3, \ldots \geq 0 \}$, with the convention that if $0 \in S$, we include only one copy of $k[1]$.
From (\ref{HH_1_case_3a}) 
$$HH_1^\s (\A) = 
( \Sigma_{s \in S'}^\oplus \; k[( b^s ,b)] ) \oplus
( \Sigma_{t \in S'}^\oplus \; k[( c^t ,c)] ) \oplus$$
$$( \Sigma_{0 \leq i \leq N}^\oplus \; k[( b^i c^{N+1-i} ,b )] ) \oplus
( \Sigma_{0 \leq i \leq  N}^\oplus \; k[( b^{i+1} c^{N-i} ,c )] ) \oplus
k [ \omega_1 ]$$
Here $S' = \{ \integers \geq N \} \cup \{ N-2, N-4, \ldots \geq 0 \}$, and
 $[ \omega_1]=[(b,c)] = - [(c,b)]$ for $N$ odd, $[\omega_1]=0$ for $N$ even.
Now, for $s$, $t \geq 0$,
$$B_0 [ b^{s+1} ] = (s+1) [ (b^s ,b)], \quad B_0 [ c^{t+1} ] = (t+1) [ (c^t ,c)]$$
Note that, for $s \geq 0$, $s+1 \in S$ if and only if $s \in S'$.
We also have 
$$B_0 [1] = [(1,1)] = [ b (1,1,1)] =0$$
By Lemma \ref{lemma:B_0_x^s_y^t}, for $0 \leq i \leq N$
 $$B_0 [ b^{i+1} c^{N+1-i} ]  = (N+1-i) [(b^{i+1} c^{N-i} , c)] + (i+1) [(b^i c^{N+1-i} ,b)]$$
Hence $\ker(B_0) = k[1]$ if $N$ is odd, 0 if $N$ is even. Further,
   \begin{equation}
 \label{HC_1_case_3}
 \label{HH_1_quotient_im_B_0}
 HH_1^\s (\A) / \im(B_0) \cong 
  \left\{
	\begin{array}{ll}
	 k^{N+1} = \Sigma_{0 \leq i \leq N}^\oplus k [ \om_1 (N,i) ]  &: N \; \even\\
	 {} & {}\\
	 k^{N+2} = \Sigma_{0 \leq i \leq N}^\oplus k [ \om_1 (N,i) ]  \oplus k[ \om_1 ] &: N \; \odd
	\end{array}\right.
\end{equation}
with generators 
\begin{equation}
\label{omega_1_N_i}
[\om_1 (N,i) ] = [(xc,b)] = [(bx,c)], \quad x= b^i c^{N-i}, \quad 0 \leq i \leq N
\end{equation}
 together with (if $N$ is odd) $[\om_1] = [(c,b)] = -[(b,c)]$.
 
 \begin{prop} 
 \label{B_1=} For $N$ odd, $B_1 =0$. For $N$ even, $\im(B_1)$ is at most one-dimensional, spanned by $[B_1 ( \om_1 (N , {\frac{1}{2}} N))]$. 
 \end{prop}
 \begin{pf}
 Recall that $HH_2^\s (\A) \cong k^{N+1}$, with generators $\om_2 (N,i)$, $0 \leq i \leq N$ given in 
 (\ref{HH_2_case_3}).
 We use the construction of $\s$-twisted Hochschild $n$-cocycles of Proposition \ref{hoch_n_cocycle}.
 Using %, (\ref{h_1_a^r+1_d^r+1}), (\ref{h_b^s+1}), 
(\ref{h_1})-(\ref{h_c^t+1}),  for each $n \in \bZ$ define a trace $h_n$ by
 \begin{equation}
\label{h_n}
h_n ( a^{i+1} ) = 0 = h_n ( d^{i+1}),  \quad h_n (b^j ) = \d_{n,j}, \quad h_n ( c^j ) = \d_{-n,j} 
\end{equation}
for $i$, $j \geq 0$. 
For the derivation $\der_b$ (\ref{der_b}) and $\s$-derivation 
$\der$ (\ref{sigma_der})
define a $\s$-twisted Hochschild 2-cocycle $\phi_{2,n}$ by
 \begin{equation}
 \label{phi_2_n}
 \phi_{2,n}  (x,y,z) = h_n ( x \der_b (y) \der(z))
 \end{equation}
 \begin{lemma} $< \phi_{2,n},   \omega_2 (N,i) > = 0$, unless $n = 2i-N$. For $i \neq {\frac{1}{2}} N$, we have $< \phi_{2,2i-N},   \omega_2 (N,i) > \neq 0$.
\end{lemma}
\begin{pf} Directly, 
$< \phi_{2,n},   \omega_2 (N,i) > = < \phi_{2,n},   q (bdx,c,a) - q^{-1} (cax,b,d)>$ (considering only potentially nonzero terms)
\begin{eqnarray}
&&=h_n ( qbdx (-c) a - q^{-1} cax b (- \l^{-1} d)) = h_n ( q^{-1} \l^{-1} caxb d - q bdxca )\nonumber\\
&&= h_n ( q^N \l^{-1} bcx ad - q^{-N} bcx da) = h_n ( q^{-N} bcx ( q^{3N+2} ad - da) )\nonumber\\
&&= q^{N+2} ( q^N - q^{-N}) h_n ( bcxad) = q^{N+2} ( q^N - q^{-N}) h_n ( bcx (1 + qbc))\nonumber\\
&&= {\frac{ q^2 ( q^{2N} - 1) ( q^2 -1)}{( q^{N+4} -1)}} h_n ( bcx)\nonumber
\end{eqnarray}
Since $bc x = b^{i+1} c^{N+1-i}$, it's clear from (\ref{h_b^s+1}), (\ref{h_c^t+1}) that $h_n (bcx) =0$ unless $n = (i+1) - (N-i+1) = 2i -N$, and $h_{2i-N} (bcx) \neq 0$ unless $2i = N$.
\end{pf}

 To find $\im(B_1)$, since $B_1 \circ B_0 =0$, we need only consider $B_1$ applied to $\om_1$ and the $\om_1 (N,i)$.
For $[\om_1] = [(c,b)] = -[(b,c)]$, which is nonzero if and only if $N$ is odd, we have
  $\deg( \om_1 ) = (0,0)$ for the $\bZ^2$-grading  (\ref{Z^2_grading}). 
  The maps $B_n$ preserve the grading, so $\deg( B_1 ( \om_1 ) ) = (0,0)$ also. 
Now, $\deg( \om_2 (N,i)) = (0,2i-N) \neq (0,0)$ for any $N$, $i$ since $N$ is odd. 
The Hochschild boundary maps $b$ also preserve the grading, so $B_1 (\om_1)$ cannot be cohomologous to any nontrivial element of $HH_2^\s (\A)$.
 
 Now consider the generators $\om_1 (N,i)$ (\ref{omega_1_N_i}).
These contain only $b$'s and $c$'s, so combining this with (\ref{B_n}), it is immediate that each $\phi_{2,n}$ vanishes on $\im (B_1)$. So for $i \neq {\frac{1}{2}}N$, we have $[ B_1 ( \om_1 (N,i))] =0$.
  \end{pf}

   The second page of the spectral sequence reads:

	\[
     	\begin{CD}
	@ V{} VV @ V{} VV @ V{} VV @ . @ . @ . @ .\\
	{0} @ <<< {HH_3^\s (\A)/ \im(B_2) } @ <<< {\ker(B_2) / \im( B_1)} @ <<< {} @ . @ . @ . @ .\\
	@ V{} VV @ V{} VV @ V{} VV @ . @ . @ . @ .\\
	{HH_3^\s (\A) / \im(B_2)  } @ <<< {\ker(B_2) / \im( B_1)} @ <<< {\ker(B_1) / \im( B_0)} @ <<< {}  @ . @ . @ . @ .\\
 	@ VVV @ V{} VV @ V{} VV @ . @ . @ . @ .\\
	{HH_2^\s (\A) / \im(B_1)} @ <<< {\ker(B_1) / \im( B_0)} @ <<< {\ker(B_0)}  @ . @ . @ . @ . @ .\\
 	@ V{} VV @ V{} VV @ . @ . @ . @ . @ .\\
{HH_1^\s (\A) / \im(B_0)} @ <<< {\ker(B_0)} @ . @ . @ . @ . @ . \\
 	@ V{} VV @ . @ .  @ . @ . @ . @ .\\
	{HH_0^\s (\A)} @ . @ . @ . @ . @ . @ .\\
    	\end{CD}
	\]\\
The only potentially nonzero differential  is 
$f : \ker(B_0) \rto HH_3^\s (\A) / \im (B_2)$, 
and after this step the spectral sequence stabilises, giving
\begin{eqnarray}
&&HC_0^\s (\A) = HH_0^\s (\A),\quad HC_1^\s (\A) = HH_1^\s (\A) / \im(B_0),\nonumber\\
&&HC_2^\s (\A) = (HH_2^\s (\A) / \im(B_1)) \oplus \ker(B_0)\nonumber\\
&&HC_{2n+3}^\s (\A) = ((HH_3^\s (\A) / \im(B_2) )/ \im f ) \oplus (\ker(B_1) / \im(B_0)),\nonumber\\
&&HC_{2n+4}^\s (\A) = (\ker(B_2) / \im( B_1)) \oplus \ker( f)\nonumber
\end{eqnarray}
 Hence $HC_0^\s (\A)$ is infinite-dimensional, given by (\ref{HH_0_case_2}), whilest $HC_1^\s (\A)$ is given by (\ref{HC_1_case_3}).
 Now, $\ker(B_0) = k[1]$, which is nonzero if and only if $N$ is odd.
  So for $N$ even, $\ker(f) = 0 = \im(f)$. For $N$ odd, $\im(B_1) =0$, hence $HC_2^\s (\A) = HH_2^\s (\A) \oplus k[1] \cong k^{N+2}$. This completes the proof of Proposition \ref{HC_case_2}.
  \end{pf}
  
Recall that $HH_2^\s(\A)$ and $HH_3^\s(\A)$ are both isomorphic to 
$k^{N+1}$ with generators  
$\omega_2 (N,i)$ (\ref{omega_2_N_i}), 
$\omega_3 (N,i)$ (\ref{omega_3_N_i}) ($0 \le i \le N$) respectively. 
By symmetry and the $\bZ^2$-grading 
($\deg(\omega_2 (N,i)) = (0, 2i-N) = \deg( \omega_3 (N,i))$) we expect,
but do not have a proof for:
\begin{conj} \begin{enumerate}
 \item{For $N$ even, $[B_1 ( \om_1 (N, {\frac{1}{2}}N))] \neq 0 \in HH_2^\s (\A)$, and is proportional to $[\om_2 (N, {\frac{1}{2}N})]$.}
 \item{For all $N$, $B_2 :  HH_2^\s (\A) / \im (B_1)  \rto HH_3^\s (\A)$ is injective. It follows that for $N$ odd, $f : \ker(B_0 ) \rto HH_3^\s (\A)$ is the zero map.}
 \end{enumerate}
 \end{conj}  
  
  From this it would follow that, for $N$ even, $HC^\s_2 (\A) \cong k^N$, generated by $[\om_2 (N,i)]$ for $i \neq {\frac{1}{2}}N$, $HC_{2n+3}^\s (\A) \cong k^{N+2}$, generated by $[\om_3 (N, {\frac{1}{2}}N)]$ together with (\ref{HH_1_quotient_im_B_0}), and $HC_{2n+4}^\s (\A) \cong k[\om_2 (N, {\frac{1}{2}}N)]$. For $N$ odd we would have $HC_{2n+3}^\s (\A) \cong k^{N+2}$, given by 
 (\ref{HH_1_quotient_im_B_0}), and $HC_{2n+4}^\s (\A) \cong k[1]$.

\subsection{Cases 3, 4 and 5}
%Old cases 6,7,8,9!

\begin{prop}
\label{cyclic_case_3} %Old case 6!
For case 3, $\mu = q^{M+1}$, $\l=q^{-(N+1)}$, $M$, $N \geq 0$,
\begin{enumerate}
\item{
$HC_0^\s(\A) \cong k^2$, with generators $[ d^{M+1} c^{N+1} ] $, $[ a^{M+1} b^{N+1} ]$.
}
\item{
$HC_1^\s(\A) \cong k^2$, with generators  $[ ( a^{M+1} b^N ,b)]$, $[ ( d^{M+1} c^N ,c)] $,\\ equivalently 
$[ ( d^M c^{N+1} ,d)] $, $[ ( a^M b^{N+1} ,a)] $.
}
\item{
$HC_n^\s(\A) =0, \quad n \geq 2$
}
\end{enumerate}
\end{prop}

\begin{pf} Recall that
$HH_0^\s (\A) \cong k^2$, generated by $[ d^{M+1} c^{N+1} ] $, $[ a^{M+1} b^{N+1} ]$,
$HH_1^\s (\A) \cong k^4$,  generated by 
$[ ( a^{M+1} b^N ,b)]$, $[ ( d^{M+1} c^N ,c)]$, $[ ( d^M c^{N+1} ,d)]$ and $[ ( a^M b^{N+1} ,a)] $,
 and $HH_2^\s (\A) \cong k^2$ with generators 
$\omega_2$ and $\omega_2^{'}$  (\ref{case6_w2}).\\

\begin{lemma}
For $B_0 : HH_0^\s (\A) \rto HH_1^\s (\A)$ we have:
\begin{eqnarray}
&&B_0 [ d^{M+1} c^{N+1} ] = (N+1) [ (d^{M+1} c^N ,c)] + (M+1) q^{- (N+1)} [ (d^M c^{N+1} ,d)]\nonumber\\
&&B_0 [ a^{M+1} b^{N+1} ] = (N+1) [ ( a^{M+1} b^N ,b)] + (M+1) q^{N+1} [ ( a^M b^{N+1} ,b)]\nonumber
\end{eqnarray}
\end{lemma}

\begin{pf} We treat only $[ d^{M+1} c^{N+1} ]$, the calculations for $[ a^{M+1} b^{N+1} ]$ are completely analogous. 
By considering $b (1, d^{M+1} c^N ,c)$, we find that
$$B_0 [d^{M+1} c^{N+1} ] = [(1, d^{M+1} c^{N+1})] = [ (d^{M+1} c^N ,c)] + q^{-(M+1)} [( c, d^{M+1} c^N )]$$
Now, for any $x$, $y \in \A$, and for all $r \geq 0$, a simple induction shows that:
\begin{equation}
\label{handy_formula}
[ (x , y^{r+1} )] = \Sigma_{j=0}^r \; [( \s( y^j ) x y^{r-j} ,y )]
\end{equation}
Also,
$[( c, d^{M+1} c^N )] = q^{M+1} [ (d^{M+1} c, c^N)] + q^{-N(M+1)} [ (c^{N+1}, d^{M+1})]$.
It follows from  (\ref{handy_formula}) that
$[ (d^{M+1} c, c^N)] = N [ (d^{M+1} c^N ,c)]$ and\\
 $[ (c^{N+1}, d^{M+1})] = (M+1) q^{M(N+1)} [ (d^M c^{N+1} ,d)]$.
Hence the result. \end{pf}

\begin{cor} It follows that:
\begin{enumerate}
\item{ $\ker(B_0) = \{ 0 \}$.}
\item{ $HH_1^\s (\A) / \im (B_0 ) \cong k^2 $, generated by $[( a^{M+1} b^N ,b)]$,  $[(d^{M+1} c^N ,c)]$,\\
equivalently $[( d^M c^{N+1} ,d)]$, $[(a^M b^{N+1} , a)]$.}
\end{enumerate}
\end{cor}

\begin{lemma}
\label{lemma:b_1_surjective}
$B_1 : HH_1^\s (\A) \rto HH_2^\s (\A)$ is surjective.
\end{lemma}
\begin{pf}
Using (\ref{B_n}) we have
\begin{eqnarray}
&B_1 ( a^M b^{N+1} , a) =& (1, a^M b^{N+1} ,a) - 
				q^{-(N+1)} (1, a, a^M b^{N+1} ) \nonumber\\
&B_1 ( d^M c^{N+1} , d) =& (1, d^M c^{N+1} ,d) 
			- q^{N+1} ( 1,d, d^M c^{N+1} ) \nonumber
\end{eqnarray} 
Consider the twisted Hochschild 2-cocycles given by 
$$\phi_2 (x,y,z) = h_{[ a^{M+1} b^{N+1}] } ( x \der_a (y) \der_b (z))$$
$$\phi_2^{'} (x,y,z) = h_{[ d^{M+1} c^{N+1}] } ( x \der_a (y) \der_b (z))$$
with $h_{[ a^{M+1} b^{N+1}] }$, $h_{[ d^{M+1} c^{N+1}] }$ defined in (\ref{case_6_HH^0}).
 Then 
 \begin{eqnarray}
&&\phi_2 ( B_1 ( a^M b^{N+1} , a) )= -(N+1) q^{-(N+1)}\nonumber\\
&&\phi_2^{'} ( B_1 ( a^M b^{N+1} , a) ) = 0 =  \phi_2 ( B_1 ( d^M c^{N+1} , d) )\nonumber\\
&&\phi_2^{'} ( B_1 ( d^M c^{N+1} , d) ) = -q^{N+1} (N+1)\nonumber
  \end{eqnarray}
  It follows that $B_1 ( a^M b^{N+1} , a)$ and $B_1 ( d^M c^{N+1} , d)$ are nontrivial and linearly independent, and hence span $HH_2^\s (\A) \cong k^2$. 
  \end{pf}

It follows that $\ker(B_1) / \im (B_0) = 0$, and $HH_2^\s (\A) / \im (B_1) =0$.
So:
\begin{enumerate}
\item{$HC_0^\s(\A) = HH_0^\s (\A)$ }
\item{$HC_1^\s(\A) = HH_1^\s (\A) / \im (B_0) \cong k^2$ }
\item{$HC_{2n+2}^\s(\A) = (HH_2^\s (\A) / \im (B_1)) \oplus \ker (B_0) = 0$ }
\item{$HC_{2n+3}^\s(\A) = \ker(B_1) / \im (B_0) = 0$}
\end{enumerate}
This completes the proof of Proposition \ref{cyclic_case_3}.
\end{pf}

Dually, we have $HC^0_\s (\A) \cong k^2$, generated by the two 0-cocycles $h_{[ a^{M+1} b^{N+1}] }$, $h_{[ d^{M+1} c^{N+1}] }$ defined in (\ref{case_6_HH^0}).
To give the generators of $HC^1_\s (\A) \cong k^2$, define a new derivation $\der' = (N+1) \der_a - (M+1) \der_b$. We have
$$\der' ( a^{M+1} b^{N+1} ) = 0 = \der' ( d^{M+1} c^{N+1} )$$
so by Lemma \ref{cyclic_1_cocycle} the twisted Hochschild 1-cocycles $\phi_1$, $\phi_1^{'}$ defined by 
$$\phi_1 (x,y) = h_{[ a^{M+1} b^{N+1}] } ( x \der' (y)), \quad 
\phi_1^{'} (x,y) = h_{[ d^{M+1} c^{N+1}] } (x \der' (y))$$
are also twisted cyclic 1-cocycles, and satisfy
$$\phi_1 ( a^M b^{N+1},a) = N+1, \quad \phi_1 ( d^M c^{N+1} ,d) = 0$$
$$\phi_1^{'} ( a^M b^{N+1},a) = 0, \quad \phi_1^{'} ( d^M c^{N+1} ,d) = -(N+1)$$
In fact $\phi_1$, $\phi_1^{'}$ are a basis for $HC^1_\s (\A) \cong k^2$.\\

\begin{prop}
\label{tch_case_4}
In case 4, %old case 8
 $\mu = q^{-(M+1)}$, $\l = q^{-(N+1)}$, $M$, $N \geq 0$ ,
\begin{enumerate}
\item{
$HC_0^\s(\A) \cong k^2$, with generators $[ d^{M+1} b^{N+1} ] $, $ [ a^{M+1} c^{N+1} ]$.
}
\item{
$HC_1^\s(\A) \cong k^2 $, generated by 
$[ ( d^{M+1} b^N ,b)] $, $[ ( a^{M+1} c^N , c) ] $,\\ equivalently 
$[ ( d^M b^{N+1} ,d)] $, $ [ ( a^M c^{N+1} ,a )]$.
}
\item{$HC_n^\s(\A) =0, \quad n \geq 2$.}
\end{enumerate}
\end{prop}

The proof is completely analogous to that of Proposition \ref{cyclic_case_3}.
We also have $HC^0_\s (\A) \cong k^2$, generated by the two 0-cocycles $h_{[ a^{M+1} c^{N+1}] }$, $h_{[ d^{M+1} b^{N+1}] }$, and $HC^1_\s (\A) \cong k^2$, generated by $\phi_1$, $\phi_1^{'}$ defined by
$$\phi_1 (x,y) = h_{[ a^{M+1} c^{N+1}] } ( x \der'' (y)), \quad 
\phi_1^{'} (x,y) = h_{[ d^{M+1} b^{N+1}] } (x \der'' (y))$$
where $\der'' = (N+1) \der_a + (M+1) \der_b$.\\

The remaining case is the trivial one:

\begin{prop} In case 5 %Old cases 7 and 9
 ($\mu = q^{\pm (M+1)}$, $M \geq 0$, $\l \notin q^{-\bN}$, and $\mu \notin q^{\bZ}$, $\l \neq 1$), we have 
$HC_n^\s(\A) = 0$ for all $n \geq 0$.
\end{prop}
\begin{pf} In each case $HH_n^\s (\A) = 0$ for all $n \geq 0$, so the spectral sequence stabilises at the first page, with all entries being zero. \end{pf}

\section{Covariant differential calculi}
\label{section_diff_calc_suq2}

In this section we identify the classes in twisted cyclic cohomology of $\A(SL_q(2))$ of the twisted cyclic cocycles arising from the three and four dimensional covariant differential calculi originally discovered for quantum $SU(2)$ by Woronowicz.

\subsection{\sc three dimensional left-covariant calculus}
\label{section:3dcalculus}

The automorphism  of $\A(SL_q (2))$ corresponding to Woronowicz's three-dimensional left-covariant calculus over quantum $SU(2)$ is
\begin{equation}
\label{3dcalc_aut}
\s(a) = q^{-2} a, \quad
\s(b) = q^4 b, \quad
\s(c) = q^{-4} c,\quad
\s(d) = q^2 d
\end{equation}
The twisted cyclic 3-cocycle over  $\A(SL_q (2))$ arising from this calculus was written down in
\cite{sw1} section 3 (denoted by $\t_{\omega,h}$) and \cite{kmt} section 5 (corresponding to the linear functional $\int$).

\begin{thm}
\label{3dcalc_cocycle_is_trivial} For $\A = \A(SL_q (2))$, we have $HC^3_\s (\A) =0$ for the automorphism (\ref{3dcalc_aut}), 
hence the twisted cyclic 3-cocycle corresponding to $\t_{\omega,h}$
and $\int$ is a trivial element of twisted cyclic cohomology. 
\end{thm}

Specializing Proposition \ref{cyclic_case_3} to the automorphism (\ref{3dcalc_aut}), we obtain:

\begin{prop}
\label{tch_3dcalc}
 For $\s = \s_{\l,\mu}$, with $\l = q^{-2}$, $\mu= q^4$, we have
\begin{enumerate}
\item{
$HC_0^\s(\A) = HH_0^\s (\A) \cong k^2$, with generators $[ d^4 c^2] $, $[ a^4 b^2]$,
}
\item{
$HC_1^\s(\A)  \cong k^2$ generated by 
$[( a^4 b,b)]$, $[(d^4 c ,c)] $,\\
equivalently $[( d^3 c^2 ,d)] $, $[(a^3 b^2 , a)]$.
}
\item{
$HC_n^\s(\A) = 0, \quad n \geq 2$
}
\end{enumerate}
\end{prop}

By duality between twisted cyclic homology and cohomology we have

\begin{cor}
$HC^0_\s (\A) \cong k^2 \cong HC^1_\s (\A)$, $HC^n_\s (\A) =0$ for $n \geq 2$. 
\end{cor}

So the twisted cyclic 3-cocycles coming from Woronowicz's three dimensional calculus 
 that appear in \cite{sw1} section 3 (denoted by $\t_{\omega,h}$) and \cite{kmt} section 5 (corresponding to the linear functional $\int$) are necessarily trivial elements of twisted cyclic cohomology, thus proving Theorem \ref{3dcalc_cocycle_is_trivial}.

\subsection{\sc four dimensional bicovariant calculi}
\label{section:4dcalculus}

It is well-known (see \cite{sw1} for example) that the twisted cyclic 4-cocycles on $\A(SL_q (2))$ coming from both Woronowicz's  four-dimensional bicovariant calculi over quantum $SU(2)$ are both simply $S^2 h$, the promotion of the twisted 0-cocycle given by the Haar functional $h$ to a 4-cocycle via the periodicity operator $S$. 
Explicitly (up to a normalising constant),
\begin{equation}
(S^2 h) (a_0 ,a_1 , a_2 ,a_3 , a_4) = h (a_0 a_1 a_2 a_3 a_4 )
\end{equation}
On basis elements, $h$ is given by
$$h( a^{i+1} b^j c^k ) = 0 = h( d^{i+1} b^j c^k )$$
\begin{equation}
\label{haar_on_pbw}
h( b^j c^k ) = \left\{
\begin{array}{cc}
(-q)^{-k} ( 1- q^{-2}) ( 1- q^{-2(k+1)} )^{-1} & : j=k\\
0 & : j \neq k 
\end{array}
\right.
\end{equation}
From (\ref{defn_mod_aut}) we see that $h$ is a well-defined $\s_{mod}^{-1}$-twisted cyclic 0-cocycle, given by $\l = q^2$, $\mu=1$, and hence corresponds to Case 1. %Was case 4!
 By inspection, we see that $h$ is exactly the twisted 0-cocycle $h_{[1]}$ defined in (\ref{h_1}).

\section{Conclusions}
\label{conclusions}

The original motivation for this work was the
belief that calculating twisted cyclic cohomology would give new insight
into existing classification results \cite{istvan,hs} for covariant
differential calculi over quantum
$SL(2)$ and quantum $SU(2)$. However, 
we see from the Woronowicz four-dimensional calculi
that  nonisomorphic calculi can give rise to cohomologous cocycles, and as  
Theorem \ref{3dcalc_cocycle_is_trivial} shows,
interesting differential calculi can correspond to trivial elements of twisted
cyclic cohomology.

The striking result that twisting via the modular automorphism overcomes the
dimension drop in Hochschild homology seems to offer the
most promising direction for future work. Similar results have been
obtained by the first author for $\Podles$ quantum spheres
\cite{tom}, and by Sitarz for
quantum hyperplanes \cite{si}. It 
seems
 natural to ask whether the
modular automorphism can overcome the dimension
drop in Hochschild homology for larger classes of quantum groups,
and look for applications of these results.  

\section{Acknowledgements}

TH: I would like to thank Gerard Murphy for originally suggesting this problem to me, and for his help and support throughout my time in Cork. 
I gratefully acknowledge the support of the EU Quantum Spaces - Noncommutative Geometry network (INP-RTN-002) and of the\\ EPSRC. 
I am also very grateful for the hospitality of the IH${\Eacute}$S 
and of the Graduiertenkolleg Quantenfeldtheorie of the Universit\"at Leipzig.\\
UK: I would like to thank M. Khalkhali, 
S. Kolb, A. Thom, C. Voigt and
especially A. Sitarz for discussions,
explanations and comments.

\end{document}